\documentclass[12pt,reqno]{amsart}

\usepackage{amscd,amssymb}
\usepackage{epsfig}
\usepackage{float}
\usepackage{url}
\usepackage[all]{xy}
\usepackage{graphicx} 
\usepackage[utf8]{luainputenc}
\newcommand{\N}{{\mathbb N}}
\newcommand{\Z}{{\mathbb Z}}
\newcommand{\Q}{{\mathbb Q}}
\newcommand{\C}{{\mathbb C}}
\newcommand{\R}{{\mathbb R}}
\renewcommand{\P}{{\mathbb P}}
\renewcommand{\H}{{\mathbb H}}

\newcommand{\CC}{{\mathcal C}}

\newcommand{\FF}{{\mathcal F}}

\newcommand{\OO}{{\mathcal O}}

\newcommand{\TT}{{\mathcal T}}

\newcommand{\www}{\widetilde}
\newcommand{\whh}{\widehat}
\newcommand{\oooo}{\overline}
\newcommand{\uuuu}{\underline}

\newcommand{\mmod}{{\rm \ mod\ }}

\renewcommand{\Re}{{\rm Re}}

\DeclareMathOperator{\Aut}{Aut}

\DeclareMathOperator{\End}{End}
\DeclareMathOperator{\id}{id}

\DeclareMathOperator{\Rad}{Rad}

\DeclareMathOperator{\Stab}{Stab}

\DeclareMathOperator{\tr}{tr}
\DeclareMathOperator{\UR}{UR}
\DeclareMathOperator{\Cl}{Cl}

\newcommand*\oline[1]{%
  \vbox{%
    \hrule height 1.2pt
    \kern0.3ex
    \hbox{%
      \ifmmode#1\else\ensuremath{#1}\fi
    }
  }
}

\begin{document}

\theoremstyle{plain}
\newtheorem{lemma}{Lemma}[section]
\newtheorem{definition/lemma}[lemma]{Definition/Lemma}
\newtheorem{theorem}[lemma]{Theorem}
\newtheorem{proposition}[lemma]{Proposition}
\newtheorem{corollary}[lemma]{Corollary}
\newtheorem{conjecture}[lemma]{Conjecture}
\newtheorem{conjectures}[lemma]{Conjectures}

\theoremstyle{definition}
\newtheorem{definition}[lemma]{Definition}
\newtheorem{withouttitle}[lemma]{}
\newtheorem{remark}[lemma]{Remark}
\newtheorem{remarks}[lemma]{Remarks}
\newtheorem{example}[lemma]{Example}
\newtheorem{examples}[lemma]{Examples}
\newtheorem{notations}[lemma]{Notations}
\newtheorem{questions}[lemma]{Questions}

\title
{Odd vanishing cycles in cyclotomic fields} 

\author{Claus Hertling and Khadija Larabi}

\address{Claus Hertling\\
Lehrstuhl f\"ur algebraische Geometrie, 
Universit\"at Mannheim,
B6 26, 68159 Mannheim, Germany}

\email{claus.hertling@uni-mannheim.de}

\address{Khadija Larabi\\
Lehrstuhl f\"ur algebraische Geometrie, 
Universit\"at Mannheim,
B6 26, 68159 Mannheim, Germany}

\email{klarabi@mail.uni-mannheim.de}

\date{January 15, 2025}





\begin{abstract}
A cusp of a Hecke group $G_q$ has
two natural lifts to the ring of integers of a cyclotomic field.
These lifts are called here odd vanishing cycles. 
All lifts of all cusps together form a discrete subset 
of $\C$ of some exquisite beauty. 
They form one or two or four orbits of a certain subgroup of the matrix 
Hecke group. The subgroup
can be considered as a monodromy group and
is an analog of a rank 2 Coxeter group, so of a dihedral group.
The paper has a research part and a larger survey part.
\end{abstract}

\maketitle

\tableofcontents

\setcounter{section}{0}

\section{Introduction}\label{s1}
\setcounter{equation}{0}
\setcounter{figure}{0}

This paper is a mixture of a survey and of a research paper.
The research part has two aims.
The first one is to present a series of discrete subsets
$\Delta^{(1)}_q$ of $\C$ for $q\in \Z_{\geq 3}$
which have some exquisite beauty. 
They come from certain natural lifts to $\C$ of the
cusps of the Hecke groups $G_q$. 
The second aim is to study an odd variant of the 
rank two Coxeter groups 
(which are of course the dihedral groups $D_{2q}$).
The odd variants are subgroups of index 1 or 2 or 4
of the matrix Hecke groups $G_q^{mat}$. 
Background of the second aim motivates us to call the
elements of a set $\Delta^{(1)}_q$ 
{\it odd vanishing cycles}. Here the {\it even vanishing cycles}
are simply the $2q$-th unit roots.

The survey part gives background material for both aims.
It consists of the sections \ref{s2}, \ref{s3}, \ref{s4},
\ref{s6}, \ref{s10} and part of section \ref{s8}.

\begin{figure}
\includegraphics[width=1.0\textwidth]{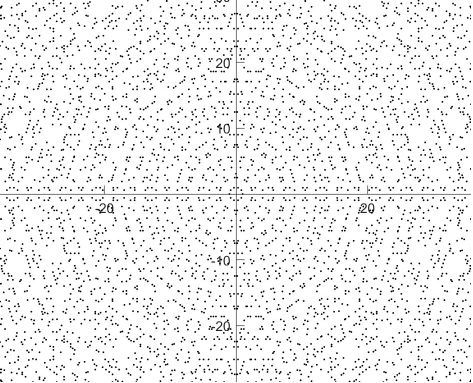}
\caption[Figure 1.1]{Part of the set $\Delta^{(1)}_5$}
\label{Fig:1.1}
\end{figure}

\begin{figure}
\includegraphics[width=1.0\textwidth]{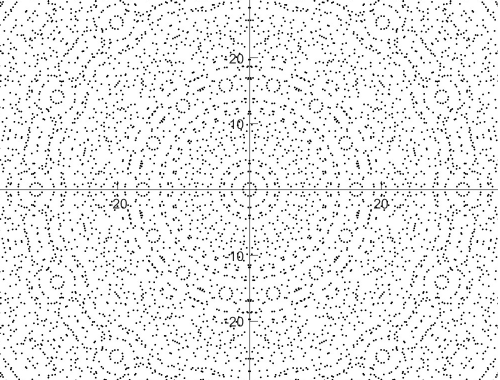}
\caption[Figure 1.2]{Part of the set $\Delta^{(1)}_7$}
\label{Fig:1.2}
\end{figure}

Figure \ref{Fig:1.1} and Figure \ref{Fig:1.2} present
the intersections of the sets $\Delta^{(1)}_5$ and 
$\Delta^{(1)}_7$ with rectangles around $0\in\C$. 

Consider a natural number $q\geq 3$, the first primitive
$2q$-th unit root $\zeta$ and two times its real part $\lambda$,
$$q\in\Z_{\geq 3},\quad\zeta:=e^{2\pi i/(2q)}\in S^1\subset\C,
\quad \lambda:=\zeta+\oooo{\zeta}=2\cos\frac{\pi}{q}\in[1,2).$$
Each two of the three matrices
$$V:=\begin{pmatrix}0&-1\\1&0\end{pmatrix},\quad
A_1:=\begin{pmatrix}1&\lambda\\0&1\end{pmatrix}
\quad\textup{and}\quad 
Q:=A_1V=\begin{pmatrix}\lambda&-1\\1&0\end{pmatrix}$$
generate the matrix group
$$G_q^{mat}:=\langle V,A_1\rangle
=\langle A_1,Q\rangle=\langle V,Q\rangle
\subset SL_2(\Z[\lambda])\subset SL_2(\R).$$
Its image $G_q$ in $PSL_2(\R)$ is a {\it Hecke group}.
The Hecke groups generalize the group $PSL_2(\Z)=G_3$.
They have been studied since Hecke's work \cite{He36}.
They are triangle groups with signature $(2,q,\infty)$
(see section \ref{s4}). 

Also the set of cusps $G_q(\infty)\subset \Q(\lambda) \,\cup\,
\{\infty\} \subset\R\,\cup\,\{\infty\}$ 
of a Hecke group has been well studied.
Section \ref{s3} will report about classical results on this set.
By definition, $G_q(\infty)$ is the image of the map
\begin{eqnarray*} 
G_q^{mat}\cdot \begin{pmatrix}1\\0\end{pmatrix}
\to \Q(\lambda)\cup\{\infty\},\ 
\begin{pmatrix}a\\c\end{pmatrix}&\mapsto& \frac{a}{c}.
\end{eqnarray*}
This map is 2:1 onto $G_q(\infty)$ because of the following. 
The structure of the triangle group $G_q$ implies that
the stabilizer of $\R\cdot(1,0)^t$  in $G_q^{mat}$ is 
$\{\pm A_1^l\,|\, l\in\Z\}$.
This easily leads to the map being 2:1 (Lemma \ref{t5.1}). 
This fact is known at least since \cite{Ro54}.

A new point of this paper is that one should embed the
column vectors of the matrices in $G_q^{mat}$ with 
the following isomorphism $\psi$ of two-dimensional 
$\R$-vector spaces into $\C$.
\begin{eqnarray*}
\psi:M_{2\times 1}(\R)\to\C,\quad 
\begin{pmatrix}x_1\\x_2\end{pmatrix}\mapsto 
x_1+\zeta^{q-1}x_2.
\end{eqnarray*}
It restricts to an isomorphism of free $\Z[\lambda]$-modules
of rank 2,
$$\psi:M_{2\times 1}(\Z[\lambda])\to\Z[\zeta].$$
The set $\Delta^{(1)}_q\subset\C$ of odd vanishing cycles is
$$\Delta^{(1)}_q:= \psi\Bigl(G_q^{mat}\cdot 
\begin{pmatrix}1\\0\end{pmatrix}\Bigr)\subset\C.$$

We denote by $\UR_{2q}$ the set of $2q$-th unit roots,
$$\UR_{2q}:=\{\zeta^k\,|\, k\in\{0,1,...,2q-1\}\}
\subset S^1\subset\C.$$
The main result of the first research aim of this paper 
is the following theorem.

\begin{theorem}\label{t1.1}
(a) The set $\Delta^{(1)}_q$ is an infinite discrete subset of 
$\C$.

(b) It is invariant under rotation by $\frac{2\pi}{2q}$.

(c) It contains the set $\UR_{2q}$ of $2q$-th unit roots.
It contains the $4q$ elements 
$\UR_{2q}\cdot \{\lambda+\zeta,1+\lambda\zeta\}$

(d) Any element of $\Delta^{(1)}_q-\UR_{2q}\cdot
\{1,\lambda+\zeta,1+\lambda\zeta\}$ has absolute value
bigger than $|\lambda+\zeta|=|1+\lambda\zeta|=\sqrt{2\lambda^2+1}$, 
with 
\begin{eqnarray*}
\begin{array}{c|c|c|c|c|c}
q & 3 & 4 & 5 & 6 & \geq 7 \\ \hline
\sqrt{2\lambda^2+1} & \sqrt{3} & \sqrt{5} & \sqrt{4+\sqrt{5}} & 
\sqrt{7} & \in]2,7375; 3[\end{array}.
\end{eqnarray*}
\end{theorem}

Our proof in section \ref{s5} is constructive. 
Following it, we produced the pictures in the 
Figures \ref{Fig:1.1}, \ref{Fig:1.2}, \ref{Fig:6.1},
\ref{Fig:6.2}, \ref{Fig:6.3}, \ref{Fig:10.1} and 
\ref{Fig:10.2}. 

The proof does not pose difficulties. It uses
the generators $A_1$ and $Q$ of $G_q^{mat}$.
Via $\psi$, the left multiplication with $Q$ on
$M_{2\times 1}(\R)$ becomes the multiplication
with $\zeta$ on $\C$. 

Theorem \ref{t1.1} should just be the first step
in studying the set $\Delta^{(1)}_q$. 
We invite the reader to go on. 
Some easy observations and some 
rather obvious naive questions are as follows: 
\begin{list}{}{}
\item[(1)]
The Figures \ref{Fig:1.1} and \ref{Fig:1.2} show a fairly
constant density of points of $\Delta^{(1)}_q$ in $\C$.
In which sense could this be made into a precise and true
statement? Lemma \ref{t6.3} says that in the cases
$q\in\{3,4,6\}$, the set $\Delta^{(1)}_q$ has 
arbitrarily large holes.
\item[(2)]
$\Delta^{(1)}_3$ is a part of a lattice $\cong\Z^2$, so
there points have a minimal distance.
Each of $\Delta^{(1)}_4$ and $\Delta^{(1)}_6$ is a part of two
overlapping lattices both $\cong\Z^2$. These sets
contain pairs of points with arbitrarily small distance,
but not triples of such points. What does hold for
the other sets $\Delta^{(1)}_q$?
\item[(3)]
The sets are not translation invariant. But the
set $\UR_{2q}$ of unit roots seems to have incomplete 
(two opposite points are missing) and blurred copies,
at least for $q=7$. Can this be made precise?
\end{list}

For the second research aim of this paper, consider
$\C=\R\cdot 1+\R\cdot i=\R\cdot 1 +\R\cdot \zeta^{q-1}$ 
as an $\R$-vector space with basis $(1,\zeta^{q-1})$, and
consider the following symmetric $\R$-bilinear form
$I^{(0)}:\C\times\C\to\R$ and the following skew-symmetric
$\R$-bilinear form $I^{(1)}:\C\times\C\to\R$,
\begin{eqnarray*}
I^{(0)}(\begin{pmatrix}1\\ \zeta^{q-1}\end{pmatrix},
(1,\zeta^{q-1}))
&=&\begin{pmatrix}2&-\lambda\\-\lambda&2\end{pmatrix},\\
I^{(1)}(\begin{pmatrix}1\\ \zeta^{q-1}\end{pmatrix},
(1,\zeta^{q-1}))
&=&\begin{pmatrix}0&-\lambda\\ \lambda&0\end{pmatrix}.
\end{eqnarray*}
$I^{(0)}$ is two times the standard scalar product with
respect to the basis $(1,i)$ of $\C$ as $\R$-vector space
(see Theorem \ref{t8.4}). 
For $a\in\C$ with $I^{(0)}(a,a)=2$ 
denote by $s^{(0)}_a$ the reflection along
$\R\cdot ia$ and by $s^{(1)}_a$ a transvection, with
\begin{eqnarray*}
s^{(k)}_a(b)=b-I^{(k)}(a,b)a.
\end{eqnarray*}
Then 
$\Gamma^{(0)}:=\langle s^{(0)}_1,s^{(0)}_{\zeta^{q-1}}\rangle$
is a Coxeter group of rank 2, namely it is isomorphic
to the dihedral group $D_{2q}$. Our odd variant is the group 
$\Gamma^{(1)}:= \langle s^{(1)}_1,s^{(1)}_{\zeta^{q-1}}\rangle$.
It turns out to be a normal subgroup of the
group $G_q^\C:=\psi\circ G_q^{mat}\circ\psi^{-1}$
of index 1 (if $q$ is odd)
or 2 (if $q\equiv 0(4)$) or 4 (if $q\equiv 2(4)$).

Section \ref{s8} gives more motivation for the dichotomy
even--odd here and for the names {\it even} and 
{\it odd vanishing cycles}. 
Section \ref{s9} gives presentations for $\Gamma^{(1)}$. 
The study of $\Gamma^{(1)}$ should be the first step towards
a theory of odd Coxeter like groups. But one rank 3 case
is discouraging (Remarks \ref{t9.5}). 

The survey part of the paper consists of the sections
\ref{s2}, \ref{s3}, \ref{s4}, \ref{s6}, \ref{s10}
and part of section \ref{s8}. 
Section \ref{s2} collects relevant known facts on cyclotomic
fields and their algebraic integers.
Section \ref{s3} presents classical results on the cusps
of the Hecke groups.
Section \ref{s4} reproves the classical fact that
$G_q$ is a triangle group with signature $(2,q,\infty)$. 
Section \ref{s6} considers the cases $q=3$, $q=4$ and $q=6$,
the only cases where $G_q$ is arithmetic.
The first part of section \ref{s8} gives a general 
construction of even and odd data, which is familiar in the
theory of isolated hypersurface singularities, and which
leads to even and odd monodromy groups and even and odd
vanishing cycles.
Section \ref{s10} discusses Rosen's \cite{Ro54}
$\lambda$-continued fractions, interesting open questions
on the set of cusps in $\Q(\lambda)\cup\{\infty\}$ and on other
$G_q$ orbits, and one result in \cite{Mc22}.

\section{Some known facts on cyclotomic fields and their
algebraic integers}\label{s2}
\setcounter{equation}{0}
\setcounter{figure}{0}

\begin{notations}\label{t2.1}
Throughout the paper we fix the following three numbers:
a natural number $q\geq 3$, the first primitive $2q$-th 
unit root $\zeta$ and two times its real part $\lambda$, 
\begin{eqnarray*}
q\in\Z_{\geq 3},\quad \zeta:=e^{2\pi i /(2q)}\in S^1\subset\C,
\quad 
\lambda:=\zeta+\oooo{\zeta}=2\cos\frac{\pi}{q}\in[1,2[ .
\end{eqnarray*}
For later use denote by $\Z[\zeta]^*$ and $\Z[\lambda]^*$ the
groups of units in $\Z[\zeta]$ and $\Z[\lambda]$, and denote
\begin{eqnarray*}
\Z[\lambda]_{>0}:=\Z[\lambda]\cap\R_{>0},\ 
\Z[\lambda]^*_{>0}:=\Z[\lambda]^*\cap \R_{>0},\ 
\Q(\lambda)_{>0}:=\Q(\lambda)\cap \R_{>0}.
\end{eqnarray*}
The cyclotomic polynomial $\Phi_{2q}\in\Z[t]$ is the minimal polynomial of $\zeta$.
It is unitary and of degree $\varphi(2q)$.
It allows to calculate the minimal polynomial $p_{min,\lambda}\in\Z[t]$ of $\lambda$, 
which is unitary and of
degree $\frac{\varphi(2q)}{2}$. In this paper $\N=\{1,2,...\}$. 
\end{notations}

\begin{examples}\label{t2.2}
\begin{eqnarray*}
\begin{array}{c|c|c|c}
q&\lambda&\Phi_{2q}&p_{min,\lambda:}\\ \hline
3 & 1 & t^2-t+1 & t-1 \\
4 & \sqrt{2} & t^4+1 & t^2-2 \\
5 & \frac{1}{2}(1+\sqrt{5}) & \sum_{j=0}^4(-t)^j & t^2-t-1 \\
6 & \sqrt{3} & t^4-t^2+1 & t^2-3 \\
7 & \approx 1,8019 & \sum_{j=0}^6(-t)^j & t^3-t^2-2t+1 
\end{array}
\end{eqnarray*}
\end{examples}

The following statements are well known. They can be found for 
example in \cite{Wa82}. 

\begin{theorem}\label{t2.3}
(a) \cite[Theorem 2.5]{Wa82}
$\Q(\zeta)$ is normal over $\Q$. Its Galois group is
$\textup{Gal}(\Q(\zeta)/\Q)\cong (\Z/2q\Z)^*$.

(b) \cite[Theorem 2.6]{Wa82} 
The ring of algebraic integers in $\Q(\zeta)$ is $\Z[\zeta]$. 

(c) \cite[Proposition 2.16]{Wa82}
The maximal real subfield of $\Q(\zeta)$ is $\Q(\lambda)$.
The ring of its algebraic integers is $\Z[\lambda]$. 

(d) \cite[Exercise 2.3, Lemma 1.6]{Wa82}
\begin{eqnarray*}
\textup{UR}_{2q}&:=& \{\zeta^k\,|\, k\in\Z\}
=\{1,\zeta,\zeta^2,...,\zeta^{2q-1}\}\\
&=& \{a\in\Q(\zeta)\,|\, \exists\ k\in\N\textup{ with }a^k=1\}\\
&=& \{a\in\Z[\zeta]\,|\, |a|=1\}
\end{eqnarray*}

(e) \cite[Proposition 2.8]{Wa82} 
If $q\notin\{2^k\,|\, k\in\N\}$, then $\zeta-1\in\Z[\zeta]^*$. 
If $q\notin\{p^k\,|\, p\textup{ a prime number},
k\in\N\}$, then $\zeta+1\in\Z[\zeta]^*$. 

(f) (Dirichlet's unit theorem, e.g. \cite[Satz 2.8.1]{Ko97} ) 
$\Z[\lambda]^*$ is isomorphic to the direct
product of $\{\pm 1\}$ and a group which is isomorphic to
$(\Z^{\varphi(2q)/2-1},+)$. 

(g) \cite[Theorem 4.12 and Corollary 4.13]{Wa82}
If $q\in\{p^k\,|\, p\textup{ a prime number},k\in\N\}$ then 
\begin{eqnarray*}
\Z[\zeta]^*&=& \textup{UR}_{2q}\cdot \Z[\lambda]^*.
\end{eqnarray*}
If $q\notin\{p^k\,|\, p\textup{ a prime number},k\in\N\}$ then
\begin{eqnarray*}
\Z[\zeta]^*&=& \textup{UR}_{2q}\cdot \{1,\zeta+1\}\cdot\Z[\lambda]^*.
\end{eqnarray*}

(h) \cite[Proposition 2.8]{Wa82} 
$\lambda$ is a unit in $\Z[\lambda]$ if and only if
$q$ is odd or $q$ is even, but 
$q\notin\{2p^k\,|\, p\textup{ prime number},\ k\in\N\}$.
\end{theorem}

Recall that two numbers $a,b\in\Z[\zeta]$ are {\it associated} 
if a unit $u\in\Z[\zeta]^*$ with $b=au$ exists. Notation:
$a\sim_{ass}b$. This is an equivalence relation.
Theorem \ref{t2.3} implies the following.

\begin{corollary}\label{t2.4}
(a) If $a$ and $b\in\Z[\zeta]$ are associated then 
$\arg(a)-\arg(b)\equiv 0\mmod \frac{2\pi}{2q}$ if
$q\in\{p^k\,|\, p\textup{ a prime number},k\in\N\}$,
and $\arg(a)-\arg(b)\equiv 0\mmod \frac{2\pi}{4q}$ if
$q\notin\{p^k\,|\, p\textup{ a prime number},k\in\N\}$.

(b) For odd $q$ one has $\Q(\lambda)=\Q(\lambda^2)
=\lambda\Q(\lambda^2)$. For even $q$ one has
$[\Q(\lambda):\Q(\lambda^2)]=2$ and 
$\Q(\lambda)\supsetneqq \lambda\Q(\lambda^2)$.

(c) $\Z[\lambda^2]$ is the ring of integers in $\Q(\lambda^2)$.

(d) \cite[Hilfssatz 2]{Le67}
For $q\in\{2p^k\,|\, p\textup{ prime number, }k\in\N\}$,
$$\lambda^{\varphi(q)}\sim_{ass} p,$$ 
i.e. $\lambda^{\varphi(q)}$ and $p$ differ only by a unit in 
$\Z[\lambda^2]$.
\end{corollary}

\begin{examples}\label{t2.5}
If $q=3$ then $\Z[\lambda]=\Z$. If $q\geq 7$ then the free part of 
$\Z[\lambda]^*$ has rank $\frac{\varphi(2q)}{2}-1\geq 2$.
If $q\in\{4,5,6\}$ then $\Z[\lambda]^*=\{\pm 1\}\cdot
\{\varepsilon_0^k\,|\, k\in\Z\}$ with 
\begin{eqnarray*}
\begin{array}{c|c|c|c}
q & 4 & 5 & 6 \\ \hline
\varepsilon_0 & \sqrt{2}+1=\lambda+1 & \frac{1}{2}(1+\sqrt{5})=\lambda &
\sqrt{3}+2=\lambda+2 
\end{array}
\end{eqnarray*}
This follows for example by applying \cite[Satz 9.5.2]{Ko97} 
to these cases. 
\end{examples}

\begin{remarks}\label{t2.6}
The norms of the field extensions $\Q(\zeta)/\Q$, $\Q(\zeta)/\Q(\lambda)$
and $\Q(\lambda)/\Q$ will be useful in section \ref{s7}. 
We denote them by 
\begin{eqnarray*}
N_{\zeta:1}:\Q(\zeta)\to\Q, &&a\mapsto
\prod_{\gamma\in\textup{Gal}(\Q(\zeta)/\Q)}\gamma(a),\\
N_{\zeta:\lambda}:\Q(\zeta)\to\Q(\lambda), &&a\mapsto a\cdot\oooo{a}=|a|^2,\\
N_{\lambda:1}:\Q(\lambda)\to\Q, &&a\mapsto 
\prod_{\gamma\in\textup{Gal}(\Q(\lambda)/\Q)}\gamma(a).
\end{eqnarray*}
Of course $N_{\zeta:1}(\Z[\zeta])\subset\Z$ and 
$N_{\zeta:1}(a)=\pm 1$ for $a\in\Z[\zeta]$ if and only if
$a\in\Z[\zeta]^*$. Similarly $N_{\lambda:1}(\Z[\lambda])\subset\Z$ and
$N_{\lambda:1}(a)=\pm 1$ for $a\in\Z[\lambda]$ if and only if 
$a\in \Z[\lambda]^*$. Finally for $a\in\Q(\zeta)$ 
\begin{eqnarray}\label{2.1}
N_{\zeta:1}(a)=N_{\lambda:1}(|a|^2).
\end{eqnarray}
\end{remarks}

$\Z[\zeta]$ respectively $\Z[\lambda]$ is a principal ideal
domain if and only if $\Q(\zeta)$ respectively $\Q(\lambda)$
has class number 1. In the case of $\Q(\zeta)$, one knows
for which $q$ this holds, in the case of $\Q(\lambda)$ not.
Theorem \ref{t2.7} and the Remarks \ref{t2.8} gives results
and remarks on the state of the art. 

\begin{theorem}\label{t2.7}
\cite[Theorem 11.1]{Wa82}
Define 
\begin{eqnarray*}
\CC_1&:=&\{1,3,5,7,9,11,13,15,17,19,21,25,27,33,35,45\}\quad\textup{and}\\
\CC_2&:=&\{1,2,3,4,5,6,7,8,9,10,11,12,15,21\}.
\end{eqnarray*}
$\Q(\zeta)$ has class number 1 if and only if 
$q\in \CC_1\cup 2\CC_2$. 
\end{theorem}

\begin{remarks}\label{t2.8}
The class number $|\Cl(\Q(\lambda))|$ is equal to 1 for
many small $q\in\Z_{\geq 3}$. van der Linden \cite{vdL82}
showed that $|\Cl(\Q(\lambda))|=1$ for each prime power $q<71$.
Schoof \cite{Sch02} restricted
$q$ to an odd prime number with $q<10000$. There are 
1285 such prime numbers. For 925 of them he expects
$|\Cl(\Q(\lambda))|=1$, including all $q<163$. \cite{vdL82} and 
\cite{Sch02} give further references. Schoof stresses that
the computation of $|\Cl(\Q(\lambda))|$ is very difficult, and that
$|\Cl(\Q(\lambda))|$ is unknown for each prime number $q\geq 71$.
\end{remarks}

\section{Classical results on the cusps of the Hecke 
groups}\label{s3}
\setcounter{equation}{0}
\setcounter{figure}{0}

The following definition and the parts (a) and (b) of
the following lemma are due to Wolfart \cite{Wo77}.
In the lemma, $|\Cl(\Q(\lambda))|$ and $|\Cl(\Q(\lambda^2))|$ 
mean the class numbers of $\Q(\lambda)$ and $\Q(\lambda^2)$.

\begin{definition}\label{t3.1} \cite{Wo77}
Let $q\in\Z_{\geq 3}$ and $\lambda$ be as above. Define
\begin{eqnarray*}
(\lambda\Q(\lambda^2))^0&:=& 
\{\frac{\lambda a}{c}\,|\, a,c\in\Z[\lambda^2],c\neq 0,(\lambda a,c)_{\Z[\lambda]}=\Z[\lambda]\}
\\
&\cup& 
\{\frac{a}{\lambda c}\,|\, a,c\in\Z[\lambda^2],c\neq 0,(a,\lambda c)_{\Z[\lambda]}=\Z[\lambda]\}.
\end{eqnarray*}
\end{definition}

\begin{lemma}\label{t3.2} 
Let $q\in\Z_{\geq 3}$ and $\lambda$ be as above.

(a) \cite{Wo77} $G_q(\infty)\subset (\lambda\Q(\lambda^2))^0\cup\{\infty\}$. 

(b) \cite{Wo77} Let $q$ be even and $(\lambda\Q(\lambda^2))^0= \lambda\Q(\lambda^2)$.
Then $|\Cl(\Q(\lambda^2))|=1$.

(c) If $PSL_2(\Z[\lambda])(\infty)=\Q(\lambda)\cup\{\infty\}$ then
$|\Cl(\Q(\lambda))|=1$.

(d) Let $q$ be odd and $G_q(\infty)=\Q(\lambda)\cup\{\infty\}$.
Then $|\Cl(\Q(\lambda))|=1$.
\end{lemma}

{\bf Proof:}
(a) By definition
$$G_q(\infty)=\{\frac{a}{c}\,|\, 
\begin{pmatrix}a&b\\c&d\end{pmatrix}\in G_q^{mat}
\textup{ for some }b,d\in\Z[\lambda]\}.$$
The group $G_q^{mat}$ is obviously a subgroup of the group
\begin{eqnarray*}
\{A=\begin{pmatrix}\lambda a&b\\c&\lambda d\end{pmatrix}\,|\, 
a,b,c,d\in\Z[\lambda^2],\det A=1\}\\
\cup\ \{A=\begin{pmatrix}a&\lambda b\\\lambda c&d\end{pmatrix}\,|\, a,b,c,d\in\Z[\lambda^2],\det A=1\}.
\end{eqnarray*}

(b) See \cite{Wo77}.

(c) The group $SL_2(\Z[\lambda])$ is a Hilbert modular group.
For such groups and in much more generality, the following
classical Lemma \ref{t3.3} holds. It implies immediately
part (c).

(d) This follows from part (c) and $G_q\subset PSL_2(\Z[\lambda])$. 
\hfill$\Box$

\begin{lemma}\label{t3.3}
(Classical, e.g. \cite[(1.1) Proposition]{vdG88})

Let $K\supset\Q$ be an algebraic number field over $\Q$.
Let $\OO\subset K$ be its ring of algebraic integers.
Let $\Cl(K)$ be its (finite) class group. 
Consider the map 
\begin{eqnarray*}
\Psi:K\cup\{\infty\}\to \Cl(K),\quad k\mapsto [(1,k)_\OO] 
\textup{ for }k\in K,\quad \infty\mapsto [(1)_\OO],
\end{eqnarray*}
which maps $k\in K$ to the class of fractional ideals 
with representative
the fractional ideal $(1,k)_\OO$ and which maps
$\infty$ to the class of $(1)_\OO=\OO$ as fractional ideal.
It is surjective. 
Its fibers are the $PSL_2(\OO)$-orbits in $K\cup\{\infty\}$.
Therefore it induces a bijective correspondence between 
the set of $PSL_2(\OO)$-orbits of elements of $K\cup\{\infty\}$
and the class group $\Cl(K)$.
\end{lemma}

The most important part of Lemma \ref{t3.2} is part (a).
The parts (b) and (d) are rather coarse because
the class numbers of $\Q(\lambda)$ and $\Q(\lambda^2)$ are equal to 1
for many small $q$ (Remark \ref{t2.8}).

The following theorem collects all known results on the
set $G_q(\infty)$ of cusps of $G_q$ for $q\in\Z_{\geq 3}$.
The first statement is Lemma \ref{t3.2} (a).
The contributions to the other statements will be explained 
in the Remarks \ref{t3.5} and in Remark \ref{t3.8} (vi).

\begin{theorem}\label{t3.4}
Let $q\in\Z_{\geq 3}$ and $\lambda$ be as above. Then
\begin{eqnarray*}
G_q(\infty)\subset (\lambda\Q(\lambda^2))^0\cup\{\infty\} &\textup{for}&
\textup{any }q\in\Z_{\geq 3}.\\
G_q(\infty)\subsetneqq (\lambda\Q(\lambda^2))^0\cup\{\infty\}&\textup{for}& 
q\in \Z_{\geq 7}-\{8,9,10,12,18,20,24\}.\\
G_q(\infty)=\lambda\Q(\lambda^2)\cup\{\infty\}&\iff&
q\in\{3,4,5,6,8,10,12\}.\\
G_q(\infty)=\Q(\lambda)\cup\{\infty\}&\iff&q\in\{3,5\}.
\end{eqnarray*}
\end{theorem}

\begin{remarks}\label{t3.5}
(i) The positive results $G_q(\infty)=\lambda\Q(\lambda^2)\cup\{\infty\}$
for $q\in\{3,4,6\}$ are classical respectively easy. 
We will discuss them in detail in section \ref{s6}.

(ii) The remarkable equality 
$G_5(\infty)=\Q(\lambda)\cup\{\infty\}$  
was first conjectured by Rosen \cite{Ro63}
(he proved in the same paper 
$G_5(\infty)\supset\{\lambda^l\,|\, l\in\Z\}$). 
Leutbecher gave a first proof of this conjecture in \cite{Le67}. 
A second proof for this case and proofs for
the positive results $G_q(\infty)=\lambda\Q(\lambda^2)\cup\{\infty\}$
for $q\in\{8,10,12\}$ were given by him in \cite{Le74}.
Very different proofs for these positive results were given
quite recently by McMullen \cite{Mc22}.
We found another proof for $q=5$, see section \ref{s7}. 

(iii) Wolfart \cite{Wo77} proved the negative results
$$G_q(\infty)\subsetneqq (\lambda\Q(\lambda^2))^0\cup\{\infty\}\textup{ for } 
q\in \Z_{\geq 7}-\{8,9,10,12,18,20,24\}.$$
He built on an idea of Borho and Rosenberger 
\cite{Bo73}\cite{BR73}, to work modulo 2. They considered
the rings $\Z[\lambda]/2\Z[\lambda]$ and $\Z[\lambda^2]/2\Z[\lambda^2]$
and the image of the group $G_q^{mat}$ in $GL_2(\Z[\lambda]/2\Z[\lambda])$,
which turns out to be a dihedral group of an order which 
divides $2q$. They considered only odd $q$ and proved
$$G_q(\infty)=\Q(\lambda)\cup\{\infty\}\Rightarrow 
q\in\{2^l+1\,|\, l\geq 1\}.$$ 
This excludes $q\in\{7,11,13,..\}$, but not $q\in\{9,17,33,..\}$.
Wolfart made more systematic use of this idea. 

(iv) Leutbecher's and Wolfart's results leave open only the 
cases $q\in\{9,18,20,24\}$.
The Remarks \ref{t3.8} will tell how these cases were solved
(negatively). 
This requires Definition \ref{t3.6} and Lemma \ref{t3.7}.
\end{remarks}

\begin{definition}\label{t3.6}
Let $q\in\Z_{\geq 3}$ and $\lambda$ be as above.
Consider a hyperbolic matrix $A\in SL_2(\Z[\lambda])$. It is called 
{\it special hyperbolic} if its fixed points are in 
$\lambda\Q(\lambda^2)$.
\end{definition}

\begin{lemma}\label{t3.7}
(Classical, e.g. 
\cite[Corollary to 2I Theorem, Ch. I 2., p 14]{Le66})
If $\Gamma\subset SL_2(\R)$ is a discrete subgroup
(i.e. a Fuchsian group), then the set of fixed points
of its parabolic elements and the set of fixed points of
its hyperbolic elements are disjoint.
\end{lemma}

\begin{remarks}\label{t3.8}
(i) In section \ref{s4} we will recover the classical
result that $G_q$ is a triangle group where the hyperbolic 
triangle is degenerate with one vertex on $\R\,\cup\,\{\infty\}$. 
This implies first that $G_q$ and $G_q^{mat}$ are discrete
groups, so Lemma \ref{t3.7} applies.
It implies second that all parabolic elements of $G_q$ are
conjugate to one another \cite[Corollary 9.2.9]{Be83}, 
so that indeed $G_q(\infty)$
is the set of fixed points of all parabolic elements of $G_q$.

(ii) $\infty$ is the fixed point of the parabolic element
$A_1$, and $0$ is the fixed point of the parabolic element
$$A_2:=\begin{pmatrix}1&0\\-\lambda&1\end{pmatrix}=VA_1V^{-1}.$$

(iii) Because of Lemma \ref{t3.7} and (ii), any hyperbolic 
element $A=\begin{pmatrix}a&b\\c&d\end{pmatrix}\in G_q^{mat}$
satisfies $bc\neq 0$. Its hyperbolic fixed points are
$$\frac{-2b}{a-d\pm\sqrt{D}}=\frac{d-a\pm\sqrt{D}}{-2c}
\in\Q(\lambda)[\sqrt{D}],$$
where $D=(\tr A)^2-4=(a+d)^2-4$ is its discriminant. 

(iv) If the group $G_q$ contains a special hyperbolic element, 
then by Lemma \ref{t3.7} 
$G_q(\infty)\subsetneqq \lambda\Q(\lambda^2)\cup\{\infty\}$. 
Therefore it is interesting
for which $q$ the group $G_q$ contains special hyperbolic 
elements. 

(v) This is not an easy problem, as the groups
$G_q^{mat}$ for $q\notin\{3,4,6\}$ are not {\it arithmetic},
i.e. not {\it commensurable} to the group $SL_2(\Z)=G_3$
(see \cite{Le67} for the notion {\it commensurable}).
For $q\notin\{3,4,6\}$, there is no known way
how to describe the matrices in $G_q^{mat}$ by equations. 
The description by the generators $V$ and $A_1$ is not 
well suited for a control of the hyperbolic elements 
and the special hyperbolic elements in $G_q^{mat}$.

(vi) Seibold \cite{Se85} gave for $q\in\{9,18,20\}$ 
special hyperbolic matrices in $G_q$ with fixed points
in $\lambda\Q(\lambda^2)$. Arnoux and Schmidt \cite{AS09} 
probably did not have \cite{Se85} and knew about it only 
from \cite[top of page 534]{Ro86} that it solved the case $q=9$ negatively. They gave special hyperbolic matrices in $G_q$ 
with fixed points in $\lambda\Q(\lambda^2)$ in the
cases $q\in\{18,20,24\}$ (and also in the case $q=14$).
The results in \cite{Se85} and \cite{AS09} together show
$G_q(\infty)\subsetneqq \lambda\Q(\lambda^2)\cup\{\infty\}$
for $q\in\{9,18,20,24\}$. 
\end{remarks}

\section{The Hecke groups are triangle groups}\label{s4}
\setcounter{equation}{0}
\setcounter{figure}{0}

This section studies the action of the Hecke group
$G_q\subset PSL_2(\Z[\lambda])$ on the upper half plane. 
It is a triangle group of signature $(2,q,\infty)$.
This is well known since Hecke's work \cite{He36}.
Theorem \ref{t4.2} establishes the triangle group
and a presentation of $G_q$. 
Corollary \ref{t4.3} gives a fundamental domain.

\begin{remarks}\label{t4.1}
(i) The images of elements of $SL_2(\R)$ and of subgroups of
$SL_2(\R)$ in $PSL_2(\R)$ will be denoted with a thick overline,
so
$$\oline{\parbox[c][0.3cm][c]{0.3cm}{.}}:
SL_2(\R)\to PSL_2(\R)$$
is the natural homomorphism, so $G_q=\oline{G_q^{mat}}$. 

(ii) Recall the following facts from the theory of 
Fuchsian groups (\cite{Be83}, especially \S 10.6).
Consider a (possibly degenerate) hyperbolic triangle 
in the upper half plane whose vertices $p_1,p_2,p_3
\in\H\cup\R\cup\{\infty\}$ are in clockwise order.
Let $l_{i,j}$ for $(i,j)\in\{(1,2),(2,3),(3,1)\}$ be the
hyperbolic line through $p_i$ and $p_j$.
Let $\alpha_i$ be the angle in the triangle at $p_i$.
If a point $p_i$ is in $\R\cup\{\infty\}$ 
then the triangle is degenerate and then $\alpha_i=0$.  

Let $r_{i,j}$ be the reflection along the circle in $\P^1\C$ 
whose part is $l_{i,j}$. It maps $\H$ to $\H$.
For $(i,j,k)\in\{(1,2,3),(2,3,1),(3,1,2)\}$, the composition
$d_j:=r_{j,k}r_{i,j}$ is in $PSL_2(\R)$. 
It is an elliptic element with fixed point $p_j$ 
and rotation angle $2\alpha_j$ if $p_j\in\H$,
and it is a parabolic element if $p_j\in\R\cup\{\infty\}$.
Obviously $d_3d_2d_1=\id$. 

If it happens that for all points $p_i$ which are in 
$\H$ (and not in $\R\cup\{\infty\}$) the angle $\alpha_i$ is 
$\alpha_i=\frac{\pi}{n_i}$
for some $n_i\in\Z_{\geq 2}$ (here $n_i:=\infty$ if $\alpha_i=0$), 
then the following holds \cite[\S 10.6]{Be83}.  
The group $\langle d_1,d_2,d_3\rangle$ acts properly
discontinuously on $\H$, and the hyperbolic triangle together
with its image under one of the three reflections
$r_{1,2},r_{2,3},r_{3,1}$ is a fundamental domain for the action
of the group on $\H$. Furthermore then the map
$d_1\mapsto\delta_1,d_2\mapsto\delta_2,d_3\mapsto\delta_3$
extends to an isomorphism from the group 
$\langle d_1,d_2,d_3\rangle$ to the group with presentation
\begin{eqnarray*}
\langle \delta_1,\delta_2,\delta_3\,|\, \delta_3\delta_2\delta_1=e,
\delta_i^{n_i}=e\textup{ for all }i\textup{ with }n_i\neq\infty
\rangle.
\end{eqnarray*}
Then the group $\langle d_1,d_2,d_3\rangle$ is called a
{\it triangle group} with signature $(n_1,n_2,n_3)$. 
\end{remarks}

The matrices $V,A_1$ and $Q\in SL_2(\Z[\lambda])$ 
were defined in the introduction, $A_2$ was defined in Remark
\ref{t3.8} (ii). We recall their definitions.
\begin{eqnarray}
V=\begin{pmatrix}0&-1\\1&0\end{pmatrix},\quad
A_1=\begin{pmatrix}1&\lambda\\0&1\end{pmatrix},\quad
Q=A_1V=\begin{pmatrix}\lambda&-1\\1&0\end{pmatrix},\label{4.1}\\
\textup{thus}\quad G_q^{mat}=\langle V,A_1\rangle
=\langle A_1,Q\rangle=\langle Q,V\rangle ,\label{4.2}\\
A_2:=VA_1V^{-1}=\begin{pmatrix}1&0\\-\lambda&1\end{pmatrix}
\textup{ with }A_1A_2=-Q^2.\label{4.3}
\end{eqnarray}

\begin{theorem}\label{t4.2}
(a) Consider the points $p_1=\infty$, $p_2=\zeta$, $p_3=i$
and $p_4=0$ and the two degenerate hyperbolic triangles
with vertices $p_1,p_2,p_3$ respectively $p_1,p_2,p_4$,
see Figure \ref{Fig:4.1}. 
The angles are denoted $\alpha_1,\alpha_2,\alpha_3$
respectively $\www{\alpha}_1,\www{\alpha}_2,\www{\alpha}_4$.
They are 
\begin{eqnarray*}
\alpha_1=\www{\alpha}_1=\www{\alpha}_4=0,
\quad \alpha_2=\frac{\pi}{q},\quad \www{\alpha}_2=\frac{2\pi}{q},
\quad \alpha_3=\frac{\pi}{2}.
\end{eqnarray*}

\begin{figure}
\includegraphics[width=0.9\textwidth]{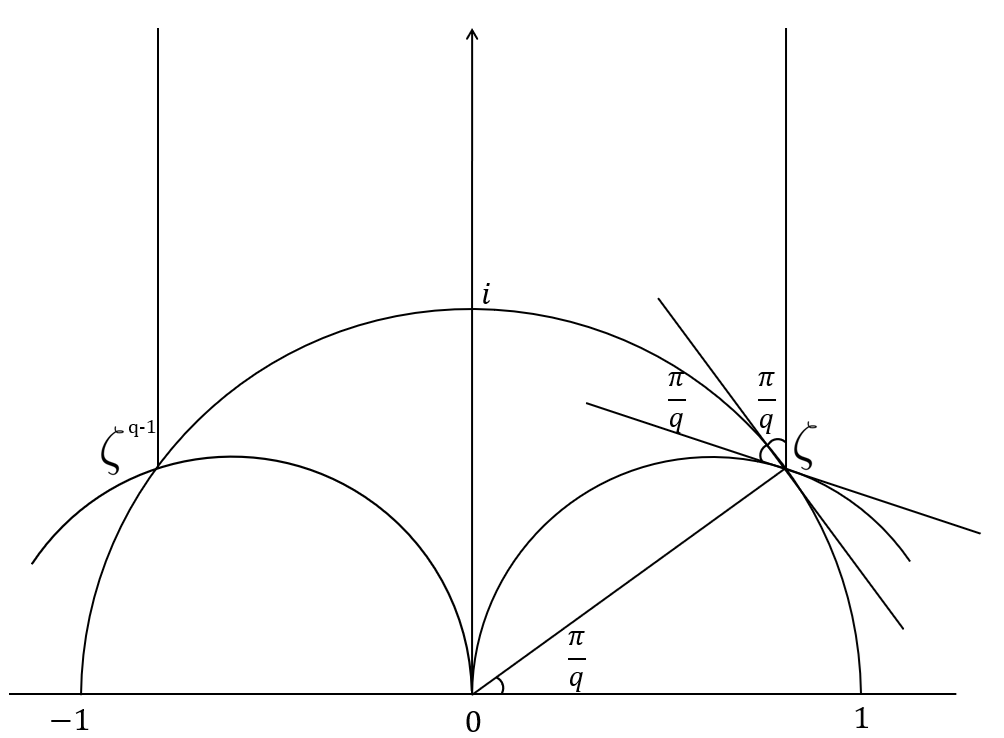}
\caption[Figure 4.1]{Some hyperbolic lines and triangles}
\label{Fig:4.1}
\end{figure}

\bigskip
(b) The elements ${\oline{V}}$, ${\oline{Q}}$,
${\oline{A_1A_2}}={\oline{Q^2}}$ are elliptic,
the elements ${\oline{A_1}},{\oline{A_2}}$ are parabolic,
their fixed points and some more information are given in the 
following table:
\begin{eqnarray*}
\begin{array}{l|c|c|c||c|c}
\textup{element} & {\oline{V}} & {\oline{Q}} &
{\oline{A_1A_2}} & {\oline{A_1}} & {\oline{A_2}} \\
\hline
\textup{fixed point} & i & \zeta & \zeta & \infty & 0 \\
\textup{rotation angle}& \pi & -\frac{2\pi}{q} &
-\frac{4\pi}{q} & - & - \\
\textup{some action} & - & - & - & \zeta^{q-1}\mapsto\zeta & 
\zeta\mapsto\zeta^{q-1} 
\end{array}
\end{eqnarray*}

(c) The group 
\begin{eqnarray}\label{4.4}
G_q=\langle {\oline{V}},{\oline{Q}}\rangle
=\langle {\oline{V}},{\oline{A_1}}\rangle
=\langle {\oline{Q}},{\oline{A_1}}\rangle
\subset PSL_2(\Z[\lambda])\subset PSL_2(\R) 
\end{eqnarray}
is the triangle group for the hyperbolic triangle 
with vertices $p_1,p_2,p_3$, so it is a triangle
group with signature $(\infty,q,2)$.
The map ${\oline{V}}\mapsto \www{v}$, 
${\oline{Q}}\mapsto \www{q}$ extends to an isomorphism
from $G_q$ to the group with presentation
\begin{eqnarray}\label{4.5}
\langle \www{v},\www{q}\,|\, 
\www{v}^2=e,\www{q}^q=e\rangle,
\end{eqnarray}
so $G_q$ is isomorphic to the free product 
$\Z_2\star\Z_q$ (with $\Z_n:=(\Z/n\Z,+)$).
\end{theorem}

{\bf Proof:}
(a) As $p_1=\infty$ and $p_4=0$, $\alpha_1=\www{\alpha}_1
=\www{\alpha}_4=0$. 
The points $\zeta^{q-1}$, $i$ and $\zeta$ are on the unit
circle, so the hyperbolic line between two of them is part of
the unit circle. The hyperbolic line between $\infty$ and $i$
is part of the verticle line through $i$. Therefore 
$\alpha_3=\frac{\pi}{2}$. 
The hyperbolic line between $\infty$ and $\zeta$ is part of
the verticle line through $\zeta$. The angle at $\zeta$ 
between this verticle line  and the unit circle is 
$\frac{\pi}{q}=\alpha_2$. See Figure \ref{Fig:4.1}. 

The map ${\oline{Q}}$ is elliptic with fixed point
$\zeta$ because
$Q\begin{pmatrix}\zeta\\ 1\end{pmatrix}=\zeta
\begin{pmatrix}\zeta\\ 1\end{pmatrix}$. 
Its rotation angle is $-2\arg\zeta=-\frac{2\pi}{q}$. 
Because $Q\begin{pmatrix}0\\ 1\end{pmatrix}=
\begin{pmatrix}-1\\ 0\end{pmatrix}$, 
${\oline{Q}}(0)=\infty$. Therefore 
${\oline{Q}}$ maps the hyperbolic line through
$0$ and $\zeta$ to the verticle line through $\zeta$.
This shows $\www{\alpha}_2=\frac{2\pi}{q}$.

(b) The statements on ${\oline{Q}}$ are proved above.
The statements on ${\oline{Q^2}}$ follow. 
The statements on ${\oline{V}}$,
${\oline{A_1}}$ and ${\oline{A_2}}$ are implied
by the following calculations.
\begin{eqnarray*}
V\begin{pmatrix}i\\ 1\end{pmatrix}
&=&i\begin{pmatrix}i\\ 1\end{pmatrix},\quad \textup{so}\quad
{\oline{V}}(i)=i\\
&&\textup{and rotation angle =}
-2\arg(i)=\pi,\\
A_1\begin{pmatrix}1\\ 0\end{pmatrix}
&=&\begin{pmatrix}1\\ 0\end{pmatrix},\quad\textup{so}\quad
{\oline{A_1}}(\infty)=\infty,\\
A_1\begin{pmatrix}\zeta^{q-1}\\ 1\end{pmatrix}
&=&\begin{pmatrix}\zeta\\ 1\end{pmatrix},\quad\textup{so}\quad
{\oline{A_1}}(\zeta^{q-1})=\zeta,
\end{eqnarray*}
\begin{eqnarray*}
A_2\begin{pmatrix}0\\ 1\end{pmatrix}
&=&\begin{pmatrix}0\\ 1\end{pmatrix},\quad\textup{so}\quad
{\oline{A_2}}(0)=0,\\
A_2\begin{pmatrix}\zeta\\ 1\end{pmatrix}
&=&(-\zeta^2)\begin{pmatrix}\zeta^{q-1}\\ 1\end{pmatrix},
\quad\textup{so}\quad {\oline{A_2}}(\zeta)=\zeta^{q-1}.
\end{eqnarray*}

(c) We use here the notations from the Remarks 
\ref{t4.1} (ii) for the vertices $p_1,p_2,p_3,p_4$ in part (a), 
the hyperbolic lines $l_{i,j}$ and the reflections $r_{i,j}$. 
Of course $r_{3,1}=r_{4,1}$. 
The parts (a) and (b) show the following.
\begin{eqnarray*}
{\oline{V}}&=& r_{3,1}r_{2,3},\quad
{\oline{Q}} = r_{1,2}r_{2,3},\\
{\oline{A_1}}&=& r_{1,2}r_{3,1},\quad
{\oline{A_2}} = r_{4,1}r_{2,4},\quad
{\oline{A_1A_2}} = r_{1,2}r_{2,4}.
\end{eqnarray*}
We recover the relation 
${\oline{Q}}={\oline{A_1}}\cdot{\oline{V}}$.
The group $G_q=\langle {\oline{V}},{\oline{Q}}\rangle
=\langle {\oline{V}},{\oline{A_1}} \rangle
=\langle {\oline{Q}},{\oline{A_1}} \rangle $
is the triangle group for the hyperbolic triangle
with vertices $p_1,p_2,p_3$. \hfill$\Box$ 

\begin{figure}
\includegraphics[width=0.95\textwidth]{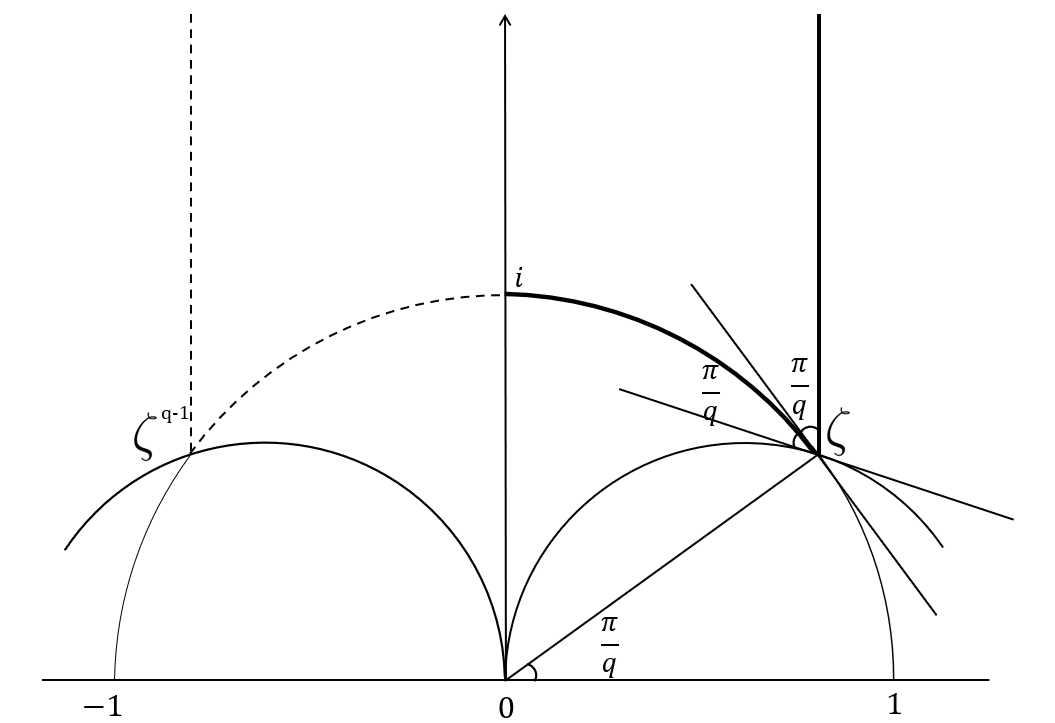}
\caption[Figure 4.2]{Fundamental domain for $G_q$}
\label{Fig:4.2}
\end{figure}

\bigskip
Corollary \ref{t4.3} states some immediate consequences
of Theorem \ref{t4.2}.
Especially, it gives a fundamental domain $\FF$ for this action. 
In Corollary \ref{t4.3} the (thin) overline over $\FF$ 
means its closure in $\C\cup\{\infty\}$. 
Figure \ref{Fig:4.2} shows this fundamental domain.

\begin{corollary}\label{t4.3}
A fundamental domain $\FF\subset\H$ for the action of $G_q$
is the set $\FF$ with 
\begin{eqnarray}\nonumber
&&\oooo{\FF}:=
\textup{the degenerate hyperbolic triangle with vertices}\\
&&\hspace*{2cm}\infty,\zeta,\zeta^{q-1},\label{4.6}\\
&&\FF\cap\{z\in\H\,|\,\Re(z)\geq 0\} 
=\oooo{\FF}\cap\{z\in\H\,|\,\Re(z)\geq 0\} ,\label{4.7}\\
&&\FF\cap \{z\in\H\,|\,\Re(z)< 0\}  
=\textup{int}(\FF)\cap \{z\in\H\,|\,\Re(z)< 0\} \label{4.8}
\end{eqnarray}
(\eqref{4.7} and \eqref{4.8} only say which points of the
boundary of $\FF$ belong to $\FF$). The stabilizers of the
points $\zeta$, $i$, $\infty$ and $0$ are
\begin{eqnarray}\label{4.9}
\Stab_{G_q}(\zeta)
=\langle{\oline{Q}}\rangle,\quad
\Stab_{G_q}(i)
=\langle{\oline{V}}\rangle,\\
\Stab_{G_q}(\infty)
=\langle{\oline{A_1}}\rangle,\quad
\Stab_{G_q}(0)
=\langle{\oline{A_2}}\rangle.\label{4.40}
\end{eqnarray}
For any $z\in\H$, the orbit $G_q\{z\}$
intersects $\FF$ in one point. 
If ${\oline{C}}\in G_q$
and $z\in\FF$ with ${\oline{C}}(z)\in\FF$ then 
${\oline{C}}(z)=z$ and 
\begin{eqnarray}\label{4.41}
\Bigl({\oline{C}}=\id\Bigr)
\textup{ or }\Bigl(z=\zeta,\ {\oline{C}}
\in \langle{\oline{Q}}\rangle\Bigr)
\textup{ or }\Bigl(z=i,\ {\oline{C}}\in 
\langle{\oline{V}}\rangle\Bigr).
\end{eqnarray}
\end{corollary}

{\bf Proof:}
Almost all statements follow immediately
from the fact that the group $G_q$ is a triangle group of the
type given in Theorem \ref{t4.2}. 

Only the statement $\Stab_{G_q}(0)
=\langle{\oline{A_2}}\rangle$ in \eqref{4.40}
requires an additional argument. It follows from the
other statement in \eqref{4.40}, from 
${\oline{V^{-1}}}(0)=\infty$ and from $A_2=VA_1V^{-1}$.
\hfill$\Box$

\section{$\Delta^{(1)}_q$ is a discrete subset of $\C$}\label{s5}
\setcounter{equation}{0}
\setcounter{figure}{0}

The matrix group $G_q^{mat}$ acts on 
$M_{2\times 1}(\Z[\lambda])$ and on $M_{2\times 1}(\R)$. 
We want actions on $\Z[\zeta]$ and $\C$. For that consider
the isomorphism of two-dimensional $\R$-vector spaces
\begin{eqnarray}\label{5.1}
\psi:M_{2\times 1}(\R)\to\C,\quad 
\begin{pmatrix}x_1\\x_2\end{pmatrix}\mapsto 
x_1+\zeta^{q-1}x_2=x_1-\zeta x_2,
\end{eqnarray}
and its restriction to an isomorphism of free
$\Z[\lambda]$-modules of rank 2,
\begin{eqnarray*}
\psi:M_{2\times 1}(\Z[\lambda])\to\Z[\zeta].
\end{eqnarray*}

\begin{lemma}\label{t5.1}
In these notations 
\begin{eqnarray}\label{5.2}
\Delta^{(1)}_q=\psi\Bigl(G_q^{mat}\cdot
\begin{pmatrix}1\\0\end{pmatrix}\Bigr).
\end{eqnarray}
The map 
\begin{eqnarray}\label{5.3}
\Delta^{(1)}_q\to G_q(\infty),\quad a+\zeta^{q-1}c\mapsto
\frac{a}{c},
\end{eqnarray}
is 2:1, so each cusp has only two lifts to $\Delta^{(1)}_q$,
and they differ only by the sign.
\end{lemma}

{\bf Proof:} The first equation is by definition. 
In order to prove that the map from $\Delta^{(1)}_q$ to
$G_q(\infty)$ is 2:1, it is sufficient to show that
the preimage of $\infty$ is $\{\pm 1\}$. 
$\Stab_{G_q}(\infty)=\langle\oline{A_1}\rangle$ in Corollary
\ref{t4.3} implies $\Stab_{G_q^{mat}}
(\R\begin{pmatrix}1\\0\end{pmatrix})
=\{\pm A_1^l\,|\, l\in\Z\}$.
The image of $\begin{pmatrix}1\\0\end{pmatrix}$ under 
$\{\pm A_1^l\,|\, l\in\Z\}$  is 
$\{\pm \begin{pmatrix}1\\0\end{pmatrix}\}$. \hfill$\Box$

\bigskip

Consider also the algebra isomorphism
\begin{eqnarray}\label{5.4}
\Psi:M_{2\times 2}(\R)&\to& \End_\R(\C),\\
B&\mapsto& \psi\circ 
(\textup{multiplication from the left with }B) \circ \psi^{-1}
\nonumber 
\end{eqnarray}
and its restriction
$$\Psi:M_{2\times 2}(\Z[\lambda])\to \End_{\Z[\lambda]}(\Z[\zeta]).$$
Denote by 
\begin{eqnarray}\label{5.5}
\mu_c:\C\to\C,\ z\mapsto cz,
\end{eqnarray} 
the multiplication with $c\in \C^*$. Observe
\begin{eqnarray}\label{5.6}
&&\Psi(Q)=\mu_\zeta,\\
\textup{namely }&&\Psi(Q):\ 1\mapsto \lambda+\zeta^{q-1}=\zeta, \ 
\zeta^{q-1}\mapsto -1=\zeta^q.\nonumber
\end{eqnarray}
Denote
\begin{eqnarray}\label{5.7}
v:=\Psi(V),\ a_1:=\Psi(A_1),\ a_2:=\Psi(A_2),\\
G_q^\C:=\Psi(G_q^{mat})=\langle v,a_1\rangle
=\langle a_1,\mu_\zeta\rangle =\langle \mu_\zeta,v\rangle
\subset\Aut_{\Z[\lambda]}(\Z[\zeta]).\label{5.8}
\end{eqnarray}
The relations $A_1V=Q=VA_2$ give the relations
\begin{eqnarray}\label{5.9}
a_1v=\mu_\zeta=va_2.
\end{eqnarray}
Lemma \ref{t5.3} describes the geometry of $a_1:\C\to\C$.
For that and also later we need notations for sectors in $\C$.

\begin{notations}\label{t5.2}
For $a,b\in S^1$ with $0<\arg(b)-\arg(a)<\pi$ 
define the open sector, the half-open sectors and the closed
sector 
\begin{eqnarray*}
S(a,b)&:=&\{x_1a+x_2b\in\C^*\,|\, x_1>0,x_2>0\},\\
S[a,b)&:=&\{x_1a+x_2b\in\C^*\,|\, x_1\geq 0,x_2>0\},\\
S(a,b]&:=&\{x_1a+x_2b\in\C^*\,|\, x_1> 0,x_2\geq 0\},\\
S[a,b]&:=&\{x_1a+x_2b\in\C^*\,|\, x_1\geq 0,x_2\geq 0\},
\end{eqnarray*}
and also the subset
\begin{eqnarray*}
S^{\geq 1}(a,b):=\{x_1a+x_2b\in\C^*\,|\, x_1\geq 1,x_2\geq 1,|x_1a+x_2b|\geq 1\}.
\end{eqnarray*}
\end{notations}

\begin{figure}
\includegraphics[width=0.7\textwidth]{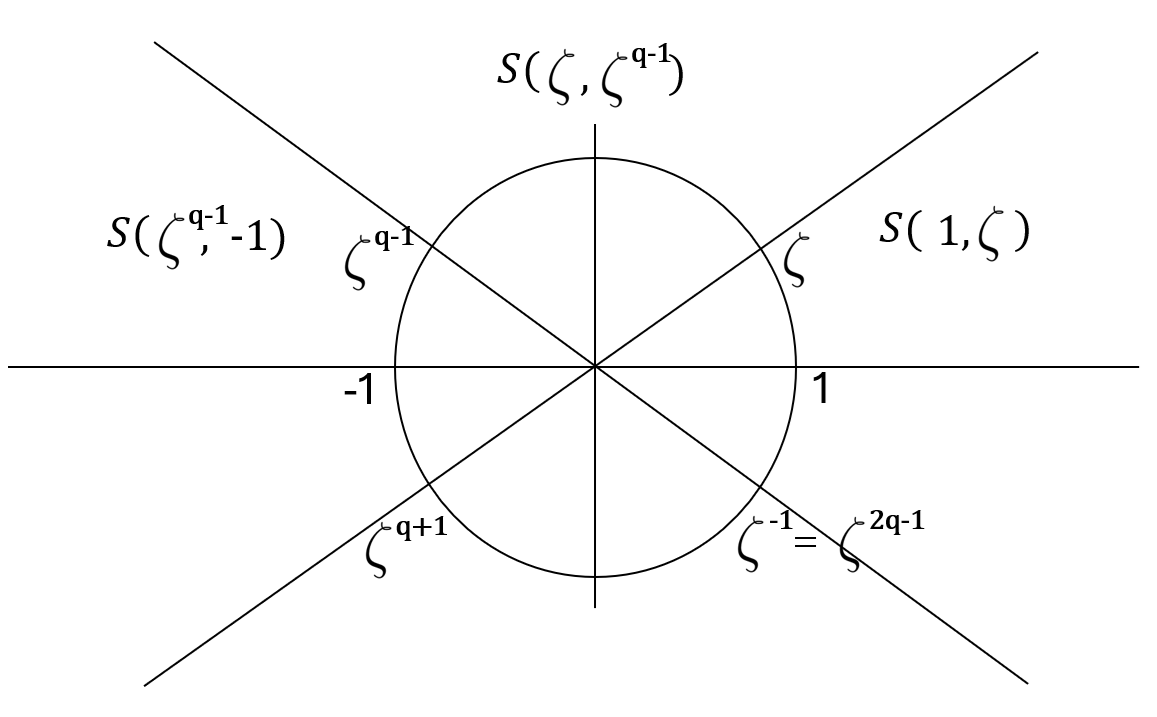}
\caption[Figure 5.1]
{The sectors $S(1,\zeta),S(\zeta,\zeta^{q-1})$
and $S(\zeta^{q-1},-1)$}
\label{Fig:5.1}
\end{figure}

\begin{figure}
\includegraphics[width=1.0\textwidth]{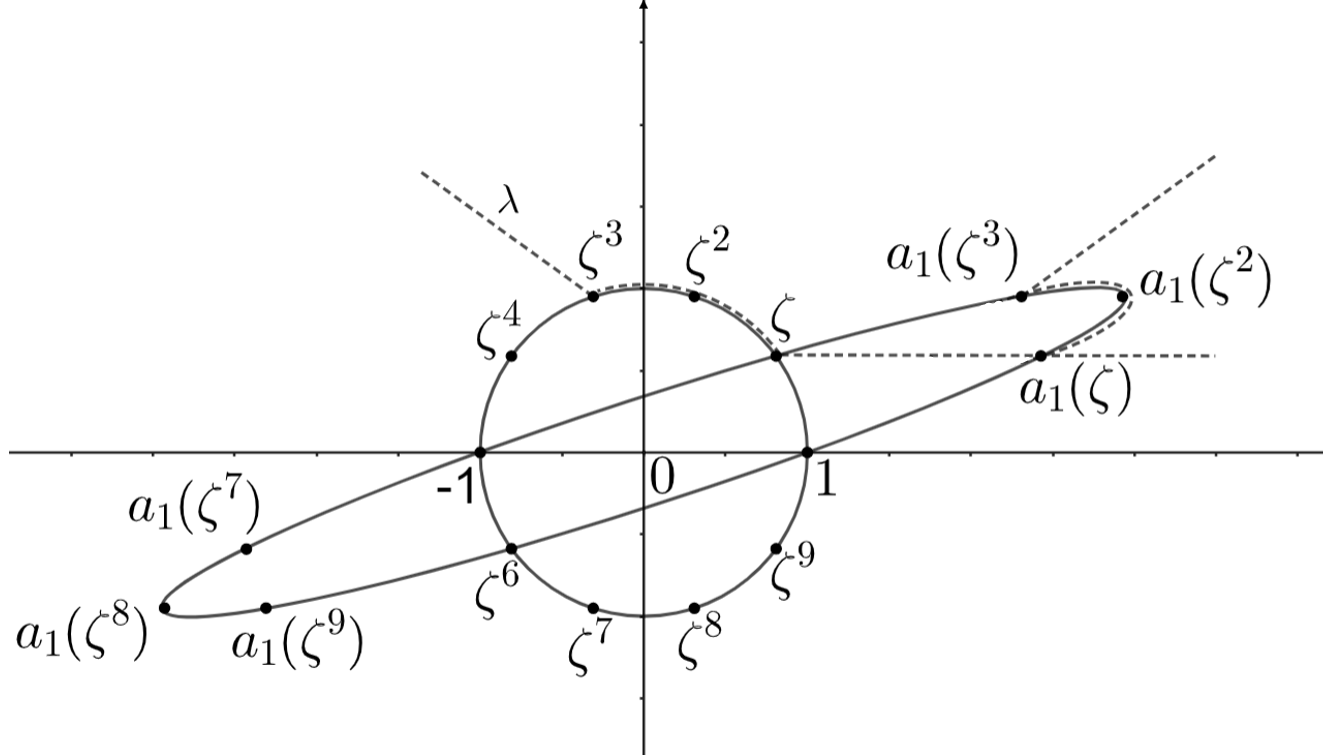}
\caption[Figure 5.2]{Action of $a_1:\C\to\C $}
\label{Fig:5.2}
\end{figure}

Figure \ref{Fig:5.1}
shows the open sectors $S(1,\zeta),S(\zeta,\zeta^{q-1})$
and $S(\zeta^{q-1},-1)$. Figure \ref{Fig:5.2} 
shows in the case $q=5$ the unit circle $S^1$, 
the 10-th unit roots and the ellipse $a_1(S^1)$.
It also shows (in form of dottes lines) the boundary of the region
$S^{\geq 1}(1,\zeta^{q-1})$ and the boundary of its image
under $a_1$, the region $a_1(S^{\geq 1}(1,\zeta^{q-1}))
\subset S(1,\zeta)$.

\bigskip
The following lemma is completely elementary, but crucial.

\begin{lemma}\label{t5.3}
The $\R$-linear automorphism $a_1:\C\to\C$ satisfies
\begin{eqnarray*}
a_1:&&1\mapsto 1,\ -1\mapsto -1,\ \zeta^{q-1}\mapsto \zeta,\\
&&S(1,\zeta^{q-1})\to S(1,\zeta)\textup{ and } 
S(\zeta^{q-1},-1)\to S(\zeta,-1),\\
|a_1(z)|&>&|z|\quad\textup{for }z\in S(1,\zeta^{q-1}),\\
|a_1(z)|&<&|z|\quad\textup{for }z\in S(\zeta^{q-1},-1).
\end{eqnarray*}
The inverse $\R$-linear automorphism $a_1^{-1}:\C\to\C$ satisfies
\begin{eqnarray*}
a_1^{-1}:&&1\mapsto 1,\ -1\mapsto -1,\ 
\zeta\mapsto \zeta^{q-1},\\
&& S(1,\zeta)\to S(1,\zeta^{q-1})\textup{ and }
S(\zeta,-1)\to S(\zeta^{q-1},-1).
\end{eqnarray*}
Additionally
\begin{eqnarray}\label{5.10}
a_1(S^{\geq 1}(1,\zeta^{q-1}))
\subset S^{\geq 1}(1,\zeta)\subset S^{\geq 1}(1,\zeta^{q-1}),\\
|a_1(z)|\geq |z|+(\sqrt{2\lambda^2+1}-1) \quad\textup{for }
z\in S^{\geq 1}(1,\zeta^{q-1}),\label{5.11}\\
\mu_\zeta^l S^{\geq 1}(1,\zeta)\subset 
S^{\geq 1}(1,\zeta^{q-1})\quad\textup{for }l\in\{1,...,q-2\}.
\label{5.12}
\end{eqnarray}
\end{lemma}

{\bf Proof:} Obviously, $a_1$ and $a_1^{-1}$ map the points 
and sectors as claimed. 

The inequality $|a_1(z)|>|z|$ for $z\in S(1,\zeta^{q-1})$ 
can be seen from Figure \ref{Fig:5.2}. 
The points $z$ in the region $S^{\geq 1}(1,\zeta^{q-1})$
with minimal difference $|a_1(z)|-|z|$ are the two vertices
$\zeta=\lambda+\zeta^{q-1}$ and 
$\zeta^{q-2}=1+\lambda\zeta^{q-1}$ 
of this region with images $a_1(\zeta)=\lambda+\zeta$ and 
$a_1(\zeta^{q-2})=1+\lambda\zeta$, 
and for them the difference is
$$|a_1(\zeta)|-|\zeta|=|a_1(\zeta^{q-2})|-|\zeta^{q-2}|=\sqrt{2\lambda^2+1}-1.$$ 
The inclusion
$a_1(S^{\geq 1}(1,\zeta^{q-1}))\subset S^{\geq 1}(1,\zeta)$
follows from
$$a_1(\alpha_1\cdot 1+\alpha_2\cdot\zeta^{q-1})=\alpha_1\cdot 1
+\alpha_2\cdot\zeta\quad\textup{for }\alpha_1\in\R_{\geq 1}
\textup{ and }\alpha_2\in\R_{\geq 1}.$$
Now $\zeta=\lambda+\zeta^{q-1}$ shows
$S^{\geq 1}(1,\zeta)\subset S^{\geq 1}(1,\zeta^{q-1})$.
Similarly \eqref{5.12} follows. The extreme case $l=q-2$
uses $\mu_\zeta^{q-2} S^{\geq 1}(1,\zeta)
=S^{\geq 1}(\zeta^{q-2},\zeta^{q-1})$ and 
$\zeta^{q-2}=1+\lambda\zeta^{q-1}$. 
\hfill$\Box$ 

\bigskip

The next lemma follows essentially from the presentation of
the Hecke group $G_q$ in \eqref{4.5} and from \eqref{4.40},
namely $\Stab_{G_q}(\infty)=\langle \oline{A_1}\rangle$,
which implies $\Stab_{G_q^\C}(1)=\langle a_1\rangle$.  
But some details in the proof are involved.

\begin{lemma}\label{t5.4}
Consider the set of tuples 
\begin{eqnarray*}
\TT(q)&:=&\Bigl\{(r,\varepsilon,
l_1,...,l_{2r})\,|\, 
r\in\N,\varepsilon\in\{0;1\},l_1\in\Z,\\
&& l_{2r}\in\{0,1,...,q-1\},\quad
\textup{in the case }r\geq 2\\\ 
&&l_2,l_4,...,l_{2r-2}\in\{1,...,q-2\}\textup{ and }
l_3,l_5,...,l_{2r-1}\in\N\Bigr\}.
\end{eqnarray*}
The map
\begin{eqnarray}
\Phi:\TT(q)&\to&G_q^\C=\langle a_1,\mu_\zeta\rangle,
\nonumber\\
(r,\varepsilon,l_1,...,l_{2r})&\mapsto& 
(-\id)^{\varepsilon}\mu_\zeta^{l_{2r}}
a_1^{l_{2r-1}}...a_1^{l_3}\mu_\zeta^{l_2}a_1^{l_1},
\label{5.13}
\end{eqnarray}
is a bijection. 
Denote $\TT^0(q):=\{(r,\varepsilon,
l_1,...,l_{2r})\in\TT(q)\,|\, l_1=0\}$. The map 
\begin{eqnarray*}
\varphi:\TT^0(q)&\to&\Delta^{(1)}_q=G_q^\C\{1\},\\
(r,\varepsilon,0,l_2,...,l_{2r})&\mapsto&
\Phi(r,\varepsilon,0,l_2,...,l_{2r})(1)
\end{eqnarray*}
is a bijection. 
\end{lemma}

{\bf Proof:}
Because of the presentation of $G_q$ in \eqref{4.5} 
and $v^2=-\id$, 
any element of $G_q^\C=\langle v,\mu_\zeta\rangle = 
\langle a_1,\mu_\zeta\rangle$ can be written 
in a unique way as a product
\begin{eqnarray}\label{5.14}
(-\id)^\delta\mu_\zeta^{k_s}
v^{-1}...v^{-1}\mu_\zeta^{k_1}v^{-\delta_2}
\end{eqnarray}
with $s\in\N$, 
$\delta,\delta_2\in\{0;1\}$,
$k_s\in\{0,1,...,q-1\}$, 
in the case $s\geq 2$ $k_1,...,k_{s-1}\in\{1,...,q-1\}$. 

Observe for $l\in\N$ 
\begin{eqnarray}\label{5.15}
a_1^l=(\mu_\zeta v^{-1})^l,\quad
a_1^{-l}=(v^{-1}\mu_\zeta^{q-1})^l.
\end{eqnarray}

Consider the product in \eqref{5.14}.
If $\delta_2=1$ and $k_1>0$, one has on the right
a factor $\mu_\zeta v^{-1}=a_1$.
If $\delta_2=0$ and $k_1=q-1$ and $s\geq 2$,
one has on the right a factor $v^{-1}\mu_\zeta^{q-1}=a_1^{-1}$.

One splits off on the right in \eqref{5.14} as many
factors $\mu_\zeta v^{-1}$ or $v^{-1}\mu_\zeta^{q-1}$
as possible. One arrives at the following unique way to
write any element of $G_q^\C$ as a product:

\begin{eqnarray}\label{5.16}
(-\id)^\delta\mu_\zeta^{m_t}
v^{-1}...v^{-1}\mu_\zeta^{m_1}a_1^{l_1}
\end{eqnarray}
with $l_1\in\Z$, $t\in\N$, 
$\delta\in\{0;1\}$,
$m_t\in\{0,1,...,q-1\}$, 
in the case $t\geq 2$ $m_2,...,m_{t-1}\in\{1,...,q-1\}$,
$m_1\in\{1,...,q-2\}$. 

In the case $t=1$ one puts $(r,\varepsilon,l_1,...,l_{2r})
:=(1,\delta,l_1,m_1)$. This leads to a tuple
in $\TT(q)$ with $r=1$, and one obtains any tuple in 
$\TT(q)$ with $r=1$ in this way. 

In the case $t\geq 2$ one rewrites each factor 
$\mu_\zeta^{m_j}v^{-1}$ for $j\in\{2,...,t\}$ 
as $\mu_\zeta^{m_j-1}a_1$.
If $m_t=0$, one rewrites $\mu_\zeta^{m_t}v^{-1}$ as 
$(-\id)\mu_\zeta^{q-1}a_1$. If some $m_j=1$, then
one obtains several factors $a_1$ behind one another.
Though it is clear that this procedure leads to a product 
in \eqref{5.13} for a tuple 
$(r,\varepsilon,l_1,...,l_{2r})\in\TT(q)$ with $r\geq 2$
and $l_2=m_1$, 
and also that one obtains any tuple in $\TT(q)$ with $r\geq 2$
by this procedure.

Therefore the map
$\Phi:\TT(q)\to G_q^\C=\langle a_1,\mu_\zeta\rangle$ 
is a bijection. 

We claim now that for two tuples in $\TT(q)$
\begin{eqnarray}\label{5.17}
\Phi((r,\varepsilon,l_1,...,l_{2r}))(1)
=\Phi((\www{r},\www{\varepsilon},\www{l}_1,...,\www{l}_{2\www{r}}))(1)
\end{eqnarray}
holds if and only if 
\begin{eqnarray}\label{5.18}
(r,\varepsilon,l_2,...,l_{2r})
=(\www{r},\www{\varepsilon},
\www{l}_2,...,\www{l}_{2\www{r}})
\end{eqnarray}
holds, so they coincide with the possible exception of the 
entries $l_{1}$ and $\www{l}_{1}$. 
The {\it if} part is clear.
In order to see the {\it only if} part, 
we go back to \eqref{5.16} and compose
one such product without $a_1^{l_1}$ from the left
with the inverse of another such product 
(where we write all indices with a hat)
without  $a_1^{\whh{l}_1}$. 
Because of $\Stab_{G_q^\C}(1)=\langle a_1\rangle$ 
it is sufficient to show that the product
\begin{eqnarray}\label{5.19}
\Bigl((-\id)^{\whh{\delta}}\mu_\zeta^{\whh{m}_{\whh{t}}}v^{-1}...
v^{-1}\mu_\zeta^{\whh{m}_1}\Bigr)^{-1}
\Bigl((-\id)^\delta\mu_\zeta^{m_t}
v^{-1}...v^{-1}\mu_\zeta^{m_1}\Bigr)
\end{eqnarray}
is a power of $a_1$ only if the indices coincide, so only
if 
$$(\whh{t},\whh{\delta},\whh{m}_1,...,\whh{m}_{\whh{t}})
=(t,\delta,m_1,...,m_t).$$
If $\whh{t}\neq t$, we can suppose $\whh{t}<t$. Then
after cancelling out as many terms in \eqref{5.19} as possible,
$\mu_\zeta^{m_1}$ remains on the right of the product. 
Then in view of \eqref{5.15} and $m_1\in\{1,...,q-2\}$
as $t\geq 2$, the product in \eqref{5.19} 
cannot be a power of $a_1$. 
The same holds if $\whh{t}=t$ and if the factors in the product
do not cancel out completely. So, the indices must coincide.

This shows the {\it only if} part in the claim above.
Therefore the map $\varphi:\TT^0(q)\to\Delta^{(1)}_q$
is a bijection. 
\hfill$\Box$

\bigskip
The bijection $\varphi:\TT^0(q)\to\Delta^{(1)}_q$ allows
to construct finite parts of the set $\Delta^{(1)}_q$
in a transparent way. This will be described in the proof of 
Lemma \ref{t5.6}. The following two notions {\it age} 
and {\it generation} will be useful.

\begin{definition}\label{t5.5}
A map $(g_1,g_2):\Delta^{(1)}_q\to\Z_{\geq 0}\times\N$ 
will be defined.
$g_1(\delta)$ will be called the {\it age} of $\delta$,
$g_2(\delta)$ will be called the {\it generation} of $\delta$.

(a) Consider a point $\delta\in\Delta^{(1)}_q$. 
By Lemma \ref{t5.4} it is the image 
\begin{eqnarray*}
\delta=\varphi(r,\varepsilon,0,l_2,...,l_{2r})
\end{eqnarray*}
under $\varphi$ of a unique tuple in $\TT^0(q)$. Then
\begin{eqnarray*}
(g_1,g_2)(\delta)&:=& (\sum_{j=2}^rl_{2j-1},r).
\end{eqnarray*}

(b) For any $q$ and any $(s,r)\in\Z_{\geq 0}\times\N$ 
define the set
\begin{eqnarray*}
\Delta^{(1,s,r)}:=\{\delta\in\Delta^{(1)}_q\,|\, 
(g_1,g_2)(\delta)=(s,r)\}.
\end{eqnarray*}
of points of age $s$ and generation $r$ and the sets
\begin{eqnarray*}
\Delta^{(1,s,*)}:=\dot\bigcup_{\rho\in\N}\Delta^{(1,s,\rho)}
\quad\textup{and}\quad
\Delta^{(1,*,r)}:=\dot\bigcup_{\sigma\in\Z_{\geq 0}}
\Delta^{(1,\sigma,r)}
\end{eqnarray*}
of points of age $s$ respectively generation $r$. 
\end{definition}

{\bf Proof of Theorem \ref{t1.1}:}
The parts (a) and (d) of Theorem \ref{t1.1} follow now
from Lemma \ref{t5.6}. 
Part (b) follows from $\mu_\zeta\in G_q^\C$.
Part (c) follows from $1\in\Delta^{(1)}_q$,
\begin{eqnarray*}
\lambda+\zeta&=&a_1(\lambda+\zeta^{q-1})=a_1(\zeta)\in\Delta^{(1)}_q,\\
1+\lambda\zeta&=&a_1(1+\lambda\zeta^{q-1})=a_1(\zeta^{q-2})\in\Delta^{(1)}_q.
\end{eqnarray*}
\hfill$\Box$

\begin{lemma}\label{t5.6}
(a) $\Delta^{(1)}_q$ is the disjoint union 
\begin{eqnarray*}
\Delta^{(1)}_q=\dot\bigcup_{(s,r)\in\Z_{\geq 0}\times\N}\Delta^{(1,s,r)}
=\dot\bigcup_{s\in\Z_{\geq 0}}\Delta^{(1,s,*)}
=\dot\bigcup_{r\in\N}\Delta^{(1,*,r)}.
\end{eqnarray*}
The set $\Delta^{(1,s,*)}$ is finite. 

(b) The set $\Delta^{(1,s,r)}$ is 
$\langle\mu_\zeta\rangle$-invariant.

(c) Write $\rho^{(s)}:=\min\{|\delta|\,|\, 
\delta\in\Delta^{(1,s,*)}\}$ for $s\in\Z_{\geq 0}$. Then
\begin{eqnarray*}
&&\rho^{(0)}=1<\rho^{(1)}=\sqrt{2\lambda^2+1},\\
&&\rho^{(s)}+(\sqrt{2\lambda^2+1}-1)\leq\rho^{(s+1)}
\quad\textup{for }s\in\Z_{\geq 0},\\
&&\rho^{(s)}\geq 1+s(\sqrt{2\lambda^2+1}-1)\quad\textup{for }
s\in\Z_{\geq 0}, \\
&&\Delta^{(1,0,*)}=\Delta^{(1,0,1)}=\textup{UR}_{2q},\\
&&\{\delta\in\Delta^{(1,1,*)}\,|\, |\delta|= \rho^{(1)}\}
=\textup{UR}_{2q}\cdot\{\lambda+\zeta,1+\lambda\zeta\}\\
&&\hspace*{4cm}=\{\delta\in\Delta^{(1)}_q-\textup{UR}_{2q}
\,|\, |\delta|\textup{ minimal}\}.
\end{eqnarray*}

(d) The set $\Delta^{(1)}_q\subset\Z[\zeta]$ is discrete.
\end{lemma}

{\bf Proof:}
(a) Obviously $\Delta^{(1)}_q$ is the disjoint union of the sets 
$\Delta^{(1,s,r)}=(g_1,g_2)^{-1}(s,r)\subset\Delta^{(1)}_q$. 
The set $\Delta^{(1,s,*)}$ is finite, because only finitely
many tuples in $\TT^0(q)$ have a given age $s$.

(b) Let $\delta\in\Delta^{(1,s,r)}$. 
We have to show $(g_1,g_2)(\zeta^k\delta)=(s,r)$. 
We know
$\delta =\Psi(r,\varepsilon,0,l_2,...,l_{2r-1},l_{2r})(1)$ with 
$l_{2r}\in\{0,1,...,q-1\}$. 
Choosing a new $\www{l}_{2r}\in\{0,1,...,q-1\}$
and a new $\www\varepsilon\in\{0;1\}$, we obtain
\begin{eqnarray*}
\Psi(r,\www{\varepsilon},0,l_2,...,l_{2r-1},\www{l}_{2r})(1)
=(-1)^{\www{\varepsilon}+\varepsilon}\mu_{\zeta}^{\www{l}_{2r}-l_{2r}}\delta\\
\textup{and}\quad
(g_1,g_2)((-1)^{\www{\varepsilon}+\varepsilon}\mu_{\zeta}^{\www{l}_{2r}-l_{2r}}
\delta)=(s,r), 
\end{eqnarray*}
so running through all possible $\www{l}_{2r}$
and both possible signs $(-1)^{\www{\varepsilon}}$, 
we see that the orbit $\langle \mu_{\zeta}\rangle\{\delta\}$
consists of elements of the same age and generation as $\delta$.

(c) First we describe informally the points in the
first and second generation.

First generation, $r=1$: $(-1)^\varepsilon\mu_\zeta^{l_2}$
with $\varepsilon\in\{0;1\}$ and $l_2\in\{0,1,...,q-1\}$
is applied to $1$. This gives all $2q$-th unit roots,
so all points in $\UR_{2q}$. 

Second generation, $r=2$: First $\mu_\zeta^{l_2}$ with
$l_2\in\{1,2,...,q-2\}$ is applied to $1$. This gives the
unit roots $\zeta,\zeta,...,\zeta^{q-2}$, so the unit roots
in $S^{\geq 1}(1,\zeta^{q-1})$ (see Figure \ref{Fig:5.2}). 
Then $a_1^{l_3}$ for some $l_3\in\N$ is applied.
By Lemma \ref{t5.3} this gives a point in 
$S^{\geq 1}(1,\zeta)\subset S(1,\zeta)$ with absolute
value $\geq 1+l_3(\sqrt{2\lambda^2+1}-1)$.
Finally $(-1)^\varepsilon\mu_\zeta^{l_4}$ with
$\varepsilon\in\{0;1\}$ and $l_4\in\{0,1,...,q-1\}$
is applied. This rotates the point by a multiple of 
$\frac{2\pi}{2q}$ and does not change its absolute value.

Higher generation, $r\geq 3$: After applying
$\mu_\zeta^{l_{2j}}a_1^{l_{2j-1}}...a_1^{l_3}\mu_\zeta^{l_2}$
to $1$ for some $j<r$, one arrives at a point in 
$S^{\geq 1}(1,\zeta^{q-1})$.
Applying $a_1^{l_{2j+1}}$ with $l_{2j+1}\in\N$ 
leads by Lemma \ref{t5.3} to a point in 
$S^{\geq 1}(1,\zeta)$ and an increase of the absolute value by
at least $l_{2j+1}(\sqrt{2l^2+1}-1)$. 
If $j<r-1$, then applying $\mu_\zeta^{l_{2j+2}}$ with 
$l_{2j+2}\in\{1,2,...,q-2\}$ rotates this point again to a point
in $S^{\geq 1}(1,\zeta^{q-1})$. 
If $j=r-1$, then applying $(-1)^\varepsilon\mu_\zeta^{l_{2r}}$
with $\varepsilon\in\{0;1\}$ and $l_{2r}\in\{0,1,...,q-1\}$
rotates the point by a multiple of $\frac{2\pi}{2q}$
and does not change its absolute value.

This shows: Any point in $\Delta^{(1)}_q$ of age
$s$ has absolute value $\geq 1+s(\sqrt{2\lambda^2+1}-1)$,
because any application of $a_1$ starts from a point
in $S^{\geq 1}(1,\zeta^{q-1})$. 
It also shows that the points of age 0 coincide with
the points of generation 1 and are the $2q$-th unit roots.
The points of age 1 are a subset of the points of generation 2,
because $1=\sum_{j=2}^rl_{2j-1}$ can hold only if $r=2$.
By Lemma \ref{t5.3} within the unit roots $\zeta,\zeta^2,...,
\zeta^{q-2}$, only the unit roots $\zeta$ and $\zeta^{q-2}$
give under the application of $a_1$ points with the
minimal absolute value $\sqrt{2l^2+1}$,
namely they give the points $a_1(\zeta)=\lambda+\zeta$
and $a_1(\zeta^{q-2})=1+\lambda\zeta$. 
Therefore the only points in $\Delta^{(1)}_q$ with this 
absolute value are the points in 
$\UR_{2q}\cdot \{\lambda+\zeta,1+\lambda\zeta\}$. 
All claims in part (c) are proved. 

(d) Any disk around 0 in $\C$ contains only finitely many elements
of $\Delta^{(1)}_q$ because each set $\Delta^{(1,s,*)}$ is finite and
because $\rho^{(s)}\geq 1+s(\sqrt{2\lambda^2+1}-1)$ grows at least
linearly with $s$. Therefore the set $\Delta^{(1)}_q$ is discrete.
\hfill$\Box$

\section{The cases $q=3$, $q=4$ and $q=6$}\label{s6}
\setcounter{equation}{0}
\setcounter{figure}{0}

The three cases $q=3$, $q=4$ and $q=6$ are the only cases
where $G_q^{mat}$ is commensurable to the group
$SL_2(\Z)=G_3^{mat}$ \cite{Le67}. For $q\notin \{3,4,6\}$
there is no known way how to describe the matrices
in $G_q^{mat}$ by equations. On the contrary, for 
$q\in\{3,4,6\}$ the group $G_q^{mat}$ and also the
set $\Delta^{(1)}_q$ can be described easily.
We give the details here.

\begin{theorem}\label{t6.1}
(a) Let $q=3$. Then $\lambda=1$ and 
\begin{eqnarray}\label{6.1}
G_3^{mat}&=&SL_2(\Z),\\ 
\Delta^{(1)}_3&=&\{a+b\zeta^2\,|\, a,b\in\Z,\gcd(a,b)=1\},\nonumber\\
G_3\{\infty\}&=&\Q\cup\{\infty\}
=\Q(\lambda)\cup\{\infty\}.\nonumber
\end{eqnarray}

(b) Let $q=4$. Then $\lambda=\sqrt{2}$ and 
\begin{eqnarray}\nonumber
G_4^{mat}&=&
\{\begin{pmatrix}a&b\sqrt{2}\\c\sqrt{2}&d\end{pmatrix}
\,|\, a,b,c,d\in\Z,ad-2bc=1\}\\
&\cup& 
\{\begin{pmatrix}a\sqrt{2}&b\\c&d\sqrt{2}\end{pmatrix}
\,|\, a,b,c,d\in\Z,2ad-bc=1\}\label{6.2} \\
\Delta^{(1)}_4&=&
\{a+c\sqrt{2}\zeta^3\,|\, a,c\in\Z,\gcd(a,2c)=1\}\nonumber \\
& \cup & 
\{b\sqrt{2}+d\zeta^3\,|\, b,d\in\Z,\gcd(2b,d)=1\},\nonumber \\
G_4\{\infty\}&=&\sqrt{2}\cdot\Q\cup\{\infty\}
=\lambda\Q(\lambda^2)\cup\{\infty\}.
\nonumber
\end{eqnarray}

(c) Let $q=6$. Then $\lambda=\sqrt{3}$ and 
\begin{eqnarray}\nonumber
G_6^{mat}&=&
\{\begin{pmatrix}a&b\sqrt{3}\\c\sqrt{3}&d\end{pmatrix}
\,|\, a,b,c,d\in\Z,ad-3bc=1\}\\
&\cup& 
\{\begin{pmatrix}a\sqrt{3}&b\\c&d\sqrt{3}\end{pmatrix}
\,|\, a,b,c,d\in\Z,3ad-bc=1\}\label{6.3} \\
\Delta^{(1)}_6&=&
\{a+c\sqrt{3}\zeta^5\,|\, a,c\in\Z,\gcd(a,3c)=1\}\nonumber \\
& \cup & 
\{b\sqrt{3}+d\zeta^5\,|\, b,d\in\Z,\gcd(3b,d)=1\},\nonumber \\
G_6\{\infty\}&=&\sqrt{3}\cdot\Q\cup\{\infty\}
=\lambda\Q(\lambda^2)\cup\{\infty\}.
\nonumber
\end{eqnarray}
\end{theorem}

{\bf Proof:} We treat the cases $q=3$, $q=4$ and $q=6$
largely together. The main point is to prove the 
explicit descriptions \eqref{6.1}, \eqref{6.2} and \eqref{6.3}
of the groups $G_q^{mat}$. 
Then the statements on $\Delta^{(1)}_q$ and $G_q\{\infty\}$
are obvious. 

The right hand sides of \eqref{6.1}, \eqref{6.2} and \eqref{6.3}
are obviously invariant under multiplication and inverting,
so they are groups. They contain $V$ and $A_1$, so they
contain $G_q^{mat}$. 

In the cases of $q=4$ and $q=6$ multiplication (e.g. from the
right) with $V$ exchanges the first and the second line of
the right hand sides of \eqref{6.2} and \eqref{6.3}. 
Therefore in all three cases
$q\in\{3,4,6\}$ it is sufficient to show that any matrix 
$C=\begin{pmatrix}a&b\lambda\\c\lambda&d\end{pmatrix}$ with $1=ad-\lambda^2bc$
is in $G_q^{mat}$.

We will show that multiplication of $C$ from the left
with $A_1^{\pm 1}$ or $A_2^{\pm 1}$ reduces one of the
entries of the first column of $C$. For any $k_1,k_2\in\Z$
\begin{eqnarray*}
A_1^{k_1}\begin{pmatrix}a\\c\lambda\end{pmatrix}
=\begin{pmatrix}a+k_1c\lambda^2\\c\lambda\end{pmatrix}
=:\begin{pmatrix}\whh{a}\\c\lambda\end{pmatrix},\\
A_2^{k_2}\begin{pmatrix}a\\c\lambda\end{pmatrix}
=\begin{pmatrix}a\\(c-k_2a)\lambda\end{pmatrix}
=:\begin{pmatrix}a\\\whh{c}\lambda\end{pmatrix}.
\end{eqnarray*}
One sees:
\begin{eqnarray}
\textup{If }0<|c\lambda^2|<|2a|\textup{ a }k_1\in\Z\textup{ exists with }
|\whh{a}|\leq \frac{1}{2}|c\lambda^2|<|a|.\label{6.4}\\
\textup{If }0<|a|<|2c|\textup{ a }k_2\in\Z\textup{ exists with }
|\whh{c}|\leq \frac{1}{2}|a|<|c|.\label{6.5}
\end{eqnarray}
If $a\neq 0$ and $c\neq 0$, at least one of the inequalities
$0<|c\lambda^2|<|2a|$ or $0<|a|<|2c|$ holds, so $|\whh{a}|<|a|$
in the first case and $|\whh{c}|<|c|$ in the second case.
Iteration leads to a matrix 
$\www{C}=
\begin{pmatrix}\www{a}&\www{b}\lambda\\ \www{c}\lambda&\www{d}\end{pmatrix}$
with $\www{a}=0$ or $\www{c}=0$.

Consider first the cases $q=4$ and $q=6$. The condition 
$1=\www{a}\www{d}-\lambda^2\www{b}\www{c}$ implies 
$\www{a}\neq 0$, so then $\www{c}=0$, $\www{a}=\www{d}=\pm 1$,
$\www{C}=\www{a}A_1^{\www{a}\www{b}}\in G_q^{mat}$
and $C\in G_q^{mat}$. 

Consider now the case $q=3$. If $\www{c}=0$ we conclude
$\www{C}\in G_3^{mat}$ and $C\in G_3^{mat}$ as in the cases 
$q\in\{4,6\}$. Suppose $\www{a}=0$. 
The matrix $V\www{C}$ has the lower left entry equal to 0.
We conclude $V\www{C}=\pm A_1^l\in G_3^{mat}$ for some $l\in\Z$
and $C\in G_3^{mat}$.
\hfill$\Box$ 

\begin{figure}
\includegraphics[width=1.0\textwidth]{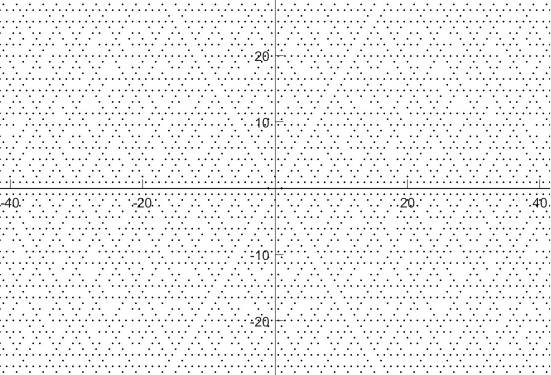}
\caption[Figure 6.1]{Part of the set $\Delta^{(1)}_3$}
\label{Fig:6.1}
\end{figure}

\begin{figure}
\vspace*{0.5cm}
\includegraphics[width=1.0\textwidth]{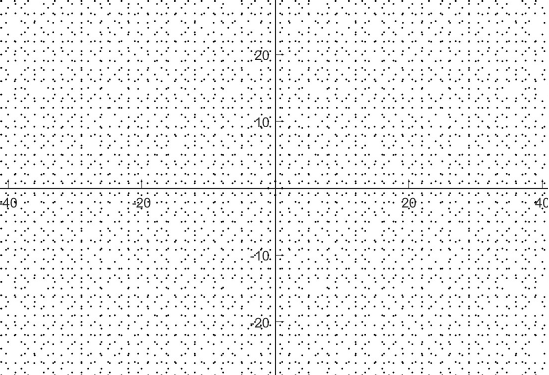}
\caption[Figure 6.2]{Part of the set $\Delta^{(1)}_4$}
\label{Fig:6.2}
\end{figure}

\begin{figure}
\vspace*{0.5cm}
\includegraphics[width=1.0\textwidth]{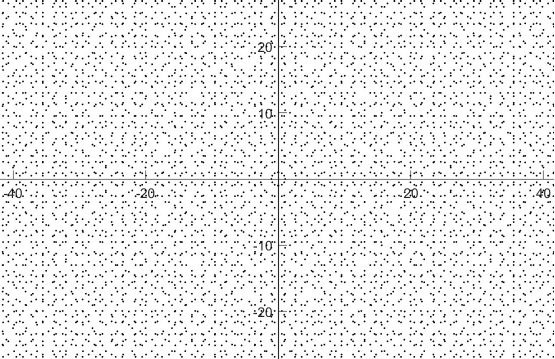}
\caption[Figure 6.3]{Part of the set $\Delta^{(1)}_6$}
\label{Fig:6.3}
\end{figure}

\begin{remarks}\label{t6.2}
(i) The Figures \ref{Fig:6.1}, \ref{Fig:6.2} and 
\ref{Fig:6.3} show in the cases
$q\in\{3,4,6\}$ the intersection of $\Delta^{(1)}_q$
with some rectangle around 0.

(ii) In the case $q=3$ $\Delta^{(1)}_q$ consists
of almost all points of the rank 2 lattice 
$\Z+\Z\zeta^2\subset\C$, namely of all points with pairwise
coprime coefficients. Especially, the minimal distance 
between different points of $\Delta^{(1)}_q$ is 1.
But this is probably the only case where this holds.

(iii) In the case $q=4$ $\Delta^{(1)}_q$ consists of almost
all points of the union of the two rank 2 lattices
$\Z+\Z\sqrt{2}\zeta^3$ and $\Z\sqrt{2}+\Z\zeta^3$.

In the case $q=6$ $\Delta^{(1)}_q$ consists of almost
all points of the union of the two rank 2 lattices
$\Z+\Z\sqrt{3}\zeta^5$ and $\Z\sqrt{3}+\Z\zeta^5$.

In both cases $\Delta^{(1)}_q$ contains pairs of points with
arbitrarily small distance between them. The reason
is that $\Z+\Z\sqrt{2}$ and $\Z+\Z\sqrt{3}$ are dense in $\R$
and contain points arbitrarily close to 0.

In both cases $\Delta^{(1)}_q$ does not contain triples of points
with arbitrarily small distances between them,
because points in one of the two lattices have a positive
minimal distance to one another.

(iv) This description of $\Delta^{(1)}_q$ as the set of almost
all points in one or two $\Z$-lattices of rank 2 lets one expect
that the density of $\Delta^{(1)}_q$ is {\it fairly constant}.
Though this is not a precise statement. The next lemma
shows that it is not obvious how to make it precise. 
\end{remarks}

\begin{lemma}\label{t6.3}
In the cases $q\in\{3,4,6\}$ there exist arbitrarily large
disks $D\subset\C$ with $D\cap\Delta^{(1)}_q=\emptyset$.
\end{lemma}

{\bf Proof:}
Choose a large number $N\in\N$. Choose $N^2$ different prime
numbers $p_{ij}$ with $i,j=1,...,N$. Define 
\begin{eqnarray*}
P_i:=\prod_{j=1}^Np_{ij},\ 
Q_j:=\prod_{i=1}^Np_{ij},\ 
\end{eqnarray*}
Then
\begin{eqnarray*}
\gcd(P_i,P_{\www{i}})
=\gcd(Q_j,Q_{\www{j}})=1\quad\textup{for }
i\neq \www{i},j\neq\www{j},\\
\gcd(P_i,Q_j)=p_{ij}.
\end{eqnarray*}
The chinese remainder theorem gives numbers $a_0,c_0\in\Z$ with
\begin{eqnarray*}
a_0&\equiv& -i\mmod P_i\quad\textup{for }i=1,...,N,\\
c_0&\equiv& -j\mmod Q_j\quad\textup{for }j=1,...,N.
\end{eqnarray*}
Both are unique up to adding any multiple of
$$P:=\prod_{i=1}^NP_i=\prod_{j=1}^NQ_j=
\prod_{i=1}^N\prod_{j=1}^Np_{ij}.$$
Because $a_0+i$ is divisible by $P_i$ and 
$c_0+j$ is divisible by $Q_j$, 
$\gcd(a_0+i,c_0+j)$ is divisible by $p_{ij}$
for $i,j\in\{1,...,N\}$. So each of the numbers
$a_0+1,a_0+2,...,a_0+N$ has a nontrivial common divisor with
each of the numbers $c_0+1,c_0+2,...,c_0+N$. 
Consider the open quadrangle with left lower vertex 0,
\begin{eqnarray*}
Q:=\left\{\begin{array}{ll}
(0,N+1)+(0,N+1)\zeta^2&\textup{ for }q=3,\\
(0,N+1)+(0,N+1)\sqrt{2}\zeta^3&\textup{ for }q=4,\\
(0,N+1)+(0,N+1)\sqrt{3}\zeta^5&\textup{ for }q=6.
\end{array}\right. 
\end{eqnarray*}

For each $l_1,l_2\in\Z$ the shifted quadrangle
\begin{eqnarray*}
\begin{array}{ll}
((a_0+l_1P)+(c_0+l_2P)\zeta^2) + Q&\textup{ for }q=3,\\
((a_0+l_1P)+(c_0+l_2P)\sqrt{2}\zeta^3) + Q&\textup{ for }q=4,\\
((a_0+l_1P)+(c_0+l_2P)\sqrt{3}\zeta^5) + Q&\textup{ for }q=6,
\end{array}
\end{eqnarray*}
does not intersect $\Delta^{(1)}_q$ in the case $q=3$,
and it does not intersect the intersection of 
$\Delta^{(1)}_q$ with the first $\Z$-lattice 
in the cases $q=4$ and $q=6$. 
This finishes the proof for the case $q=3$.
Here the first and second $\Z$-lattice are
$\Z+\Z\sqrt{2}\zeta^3$ and $\Z\sqrt{2}+\Z\zeta^3$ in the case
$q=4$, and they are 
$\Z+\Z\sqrt{3}\zeta^5$ and $\Z\sqrt{3}+\Z\zeta^5$ in the case
$q=6$. 

In the cases $q=4$ and $q=6$, 
we have to construct an overlap of one of these holes 
with a similar hole in the 
intersection of $\Delta^{(1)}_q$ with the second $\Z$-lattice. 

We can choose $N^{(2)}$ sufficiently large, e.g.
$N^{(2)}>\sqrt{3}P$, and construct
holes in the intersection of $\Delta^{(1)}_q$ with the
second $\Z$-lattice which are so large that they contain
at least one of the holes above in the intersection
of $\Delta^{(1)}_q$ with the first $\Z$-lattice.
\hfill$\Box$

\section{The case $q=5$}\label{s7}
\setcounter{equation}{0}
\setcounter{figure}{0}

In this section we offer our own proof for the beautiful result
(part of Theorem \ref{t3.4}) 
$G_q(\infty)=\Q(\lambda)\cup\{\infty\}
\textup{ in the case } q=5.$
This is the fourth proof, after two proofs by Leutbecher
\cite{Le67}\cite{Le74} and one by McMullen \cite{Mc22}.
It is not close to the other proofs.
Theorem \ref{t7.1} reformulates the statement.
The implication \eqref{7.1} 
$\Rightarrow G_5(\infty)=\Q(\lambda)\cup\{\infty\}$ is obvious.
The inverse implication uses Lemma \ref{t5.1} and that
$\Q(\lambda)$ is a principal ideal domain.

\begin{theorem}\label{t7.1}
Let $q=5$. Then
\begin{eqnarray}\label{7.1}
\Z[\zeta]=\{0\}\, \cup\, \Z[\lambda]_{>0}\cdot \Delta^{(1)}_5
=\{0\}\ \dot\cup\  
\dot\bigcup_{\delta\in\Delta^{(1)}_5}\Z[\lambda]_{>0}\cdot\delta.
\end{eqnarray}
\end{theorem}

{\bf Proof:} Recall that by Theorem \ref{t2.3} (g)
the set $\Z[\zeta]^*$ of units in $\Z[\zeta]$ is 
\begin{eqnarray*}
\Z[\zeta]^*=\Z[\lambda]^*\cdot \textup{UR}_{10}
=\Z[\lambda]^*_{>0}\cdot\textup{UR}_{10}
\end{eqnarray*} 
and that $\gamma\in\Z[\zeta]$ is a unit (so in $\Z[\zeta]^*$) 
if and only if $N_{\zeta:1}(\gamma)=\pm 1$.
Also recall from Theorem \ref{t2.3} (c)
$\Z[\zeta]\cap\R=\Z[\lambda]$, which implies 
\begin{eqnarray*}
\Z[\zeta]\cap \bigl(\R_{>0}\cdot \textup{UR}_{10}\bigr)
=\Z[\lambda]_{>0}\cdot \textup{UR}_{10}.
\end{eqnarray*}
By Lemma \ref{t5.1} a real half-line $\R_{>0}\cdot\gamma$
for $\gamma\in\Z[\zeta]-\{0\}$ contains at most one 
point in $\Delta^{(1)}_5$. 
Therefore for the proof of Theorem \ref{t7.1}
it is sufficient to prove 
$\Z[\zeta]-\{0\}=\Z[\lambda]_{>0}\cdot \Delta^{(1)}_5$.

The idea of the proof is as follows.
Let $\gamma\in\Z[\zeta]-\bigl(
\{0\}\cup\Z[\lambda]_{>0}\cdot\textup{UR}_{10}\bigr)$. 
Then $\gamma$ is not a unit, so $|N_{\zeta:1}(\gamma)|>1$.
We will describe a procedure (P) which gives a group element
$c\in G_q^\C$ such that 
\begin{eqnarray}\label{7.2}
|N_{\zeta:1}(c(\gamma))|<|N_{\zeta:1}(\gamma)|.
\end{eqnarray}
The procedure (P) can be iterated only finitely often,
because $N_{\zeta:1}$ has values in $\Z$. 
So a finite number of iterations yields a group
element $b\in G_q^\C$ such that
\begin{eqnarray*}
b(\gamma)&\in&\Z[\lambda]_{>0}\cdot\textup{UR}_{10},\\
\textup{so } b(\gamma)&=&u\cdot \zeta^k\textup{ for some }
k\in\{0,1,...,9\}\textup{ and some }u\in\Z[\lambda]_{>0},\\
\textup{so }\gamma&=& u\cdot b^{-1}(\zeta^k)
\textup{ and }b^{-1}(\zeta^k)\in\Delta^{(1)}_5.
\end{eqnarray*}
This implies $\Z[\zeta]=\Z[\lambda]_{>0}\cdot \Delta^{(1)}_5$.

Now we describe the procedure (P) and then prove that it works. 
\bigskip

{\bf Procedure (P):}
{\it Let $\gamma\in\Z[\zeta]-\bigl(
\{0\}\cup\Z[\lambda]_{>0}\cdot \textup{UR}_{10}\bigr)$. 
There is a unique $k\in\{0,1,...,9\}$ with 
\begin{eqnarray}\label{7.3}
\mu_\zeta^k(\gamma)\in 
S(1,\zeta^{0,5}]\ \dot\cup\ S(\zeta^{4,5},-1).
\end{eqnarray}
If $\mu_\zeta^k(\gamma)\in S(1,\zeta^{0,5}]$ there is a
unique $n\in\N$ with 
\begin{eqnarray}\label{7.4}
a_1^{-n}(\mu_\zeta^k(\gamma))\in S[\zeta,\zeta^4).
\end{eqnarray}
If $\mu_\zeta^k(\gamma)\in S(\zeta^{4,5},-1)$ there is a 
unique $n\in\Z_{<0}$ with 
\begin{eqnarray}\label{7.5}
a_1^{-n}(\mu_\zeta^k(\gamma))\in S(\zeta,\zeta^4].
\end{eqnarray}
In both cases define $b:=a_1^{-n}\mu_\zeta^k$. We claim
\begin{eqnarray}\label{7.6}
|N_{\zeta:1}(b(\gamma))|<|N_{\zeta:1}(\gamma)|.
\end{eqnarray}
}

In \eqref{7.4}, for $\mu_\zeta^k(\gamma)$ close to 
$\R_{>0}\zeta^{0,5}$ $n=1$, for $\mu_\zeta^k(\gamma)$
close to $\R_{>0}\cdot 1$ $n\in\N$ gets large.
Analogously in \eqref{7.5}.

\begin{figure}
\includegraphics[width=0.5\textwidth]{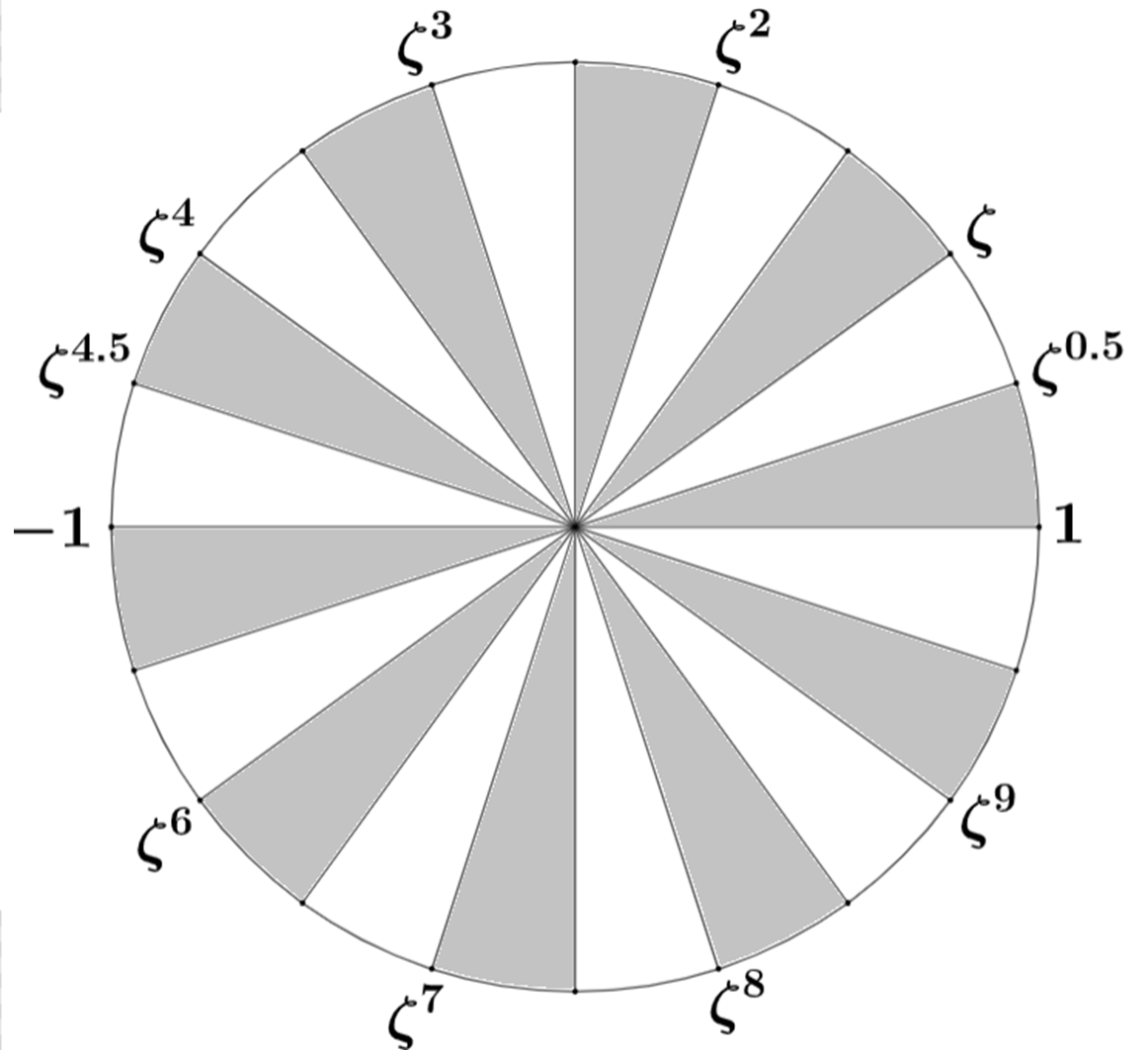}
\hspace*{0.3cm}
\includegraphics[width=1.0\textwidth]{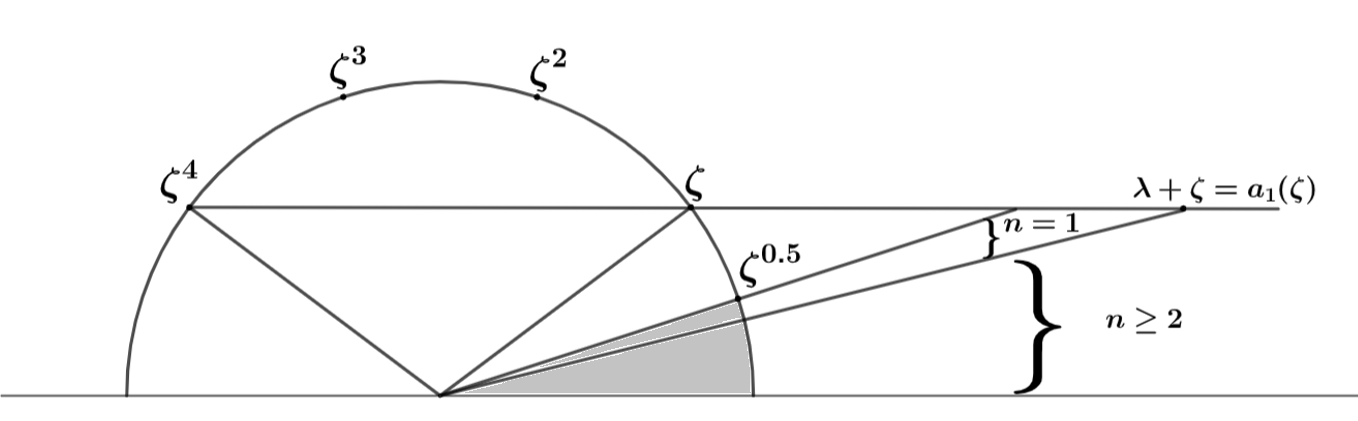}
\caption[Figure 7.1]
{20 sectors and the cases $n=1$ and $n\geq 2$ of 
Procedure (P)}\label{Fig:7.1}
\end{figure}

The two pictures in Figure \ref{Fig:7.1}
illustrate this procedure (P)
and a part of the proof of Theorem \ref{t7.1}. 

The existence and uniqueness of $k$ is easy: As 
$\gamma\notin\R_{>0}\cdot \textup{UR}_{10}$, 
there is a unique $\www{k}\in\{0,1,...,9\}$ with 
$\mu_\zeta^{\www{k}}(\gamma)\in S(1,\zeta)$.
In the case 
$\mu_\zeta^{\www{k}}(\gamma)\in S(\zeta^{0,5},\zeta)$
one moves this element with $\mu_\zeta^4$ to 
$S(\zeta^{4,5},-1)$. 

Recall that $a_1^{-1}$ maps $S(1,\zeta)$ to $S(1,\zeta^4)$
and that $a_1$ maps $S(\zeta^4,-1)$ to $S(\zeta,-1)$.
This shows existence and uniqueness of $n\in\N$ in the case
\eqref{7.4} and of $n\in\Z_{<0}$ in the case \eqref{7.5}. 

From now on we will restrict to the case \eqref{7.4}.
The case \eqref{7.5} is completely analogous
(except that the half-line $\R_{>0}\cdot \zeta^{0,5}$ is 
included, while the half-line $\R_{>0}\cdot \zeta^{4,5}$
is excluded). We observe:
\begin{eqnarray*}
n=1&\Rightarrow& 
\mu_\zeta^k(\gamma)\in 
S[\frac{\lambda+\zeta}{|\lambda+\zeta|},\zeta^{0,5}],
\ a_1^{-1}(\mu_\zeta^k(\gamma))\in S[\zeta,\zeta^2],\\
n\geq 2&\Rightarrow& 
\mu_\zeta^k(\gamma)\in 
S(1, \frac{\lambda+\zeta}{|\lambda+\zeta|}),
\ a_1^{-n}(\mu_\zeta^k(\gamma))\in S[\zeta,\zeta^4),
\end{eqnarray*}
because $\zeta=\lambda+\zeta^4$, $a_1(\zeta)=\lambda+\zeta$, 
$\zeta^2=\lambda+\lambda\zeta^4$ and 
$a_1(\zeta^2)=\lambda+\lambda\zeta=\lambda\sqrt{2+\lambda}\cdot\zeta^{0,5}$. 

{\bf The case $n=1$:} The element 
$b(\gamma)=a_1^{-1}(\mu_\zeta^k(\gamma))\in 
S[\zeta,\zeta^2]$ can be written as 
\begin{eqnarray*}
b(\gamma)=a_1^{-1}(\mu_\zeta^k(\gamma))= 
\kappa(x_1+\zeta^4)&&\textup{with suitable }
\kappa\in\Z[\lambda]_{>0}\\
&&\textup{and } x_1\in \kappa^{-1}\Z[\lambda]\cap [1,\lambda],
\end{eqnarray*}
here $x_1\in[1,\lambda]$ because 
$x_1+\zeta^4\in S[\zeta,\zeta^2]$, 
$\lambda+\zeta^4=\zeta$ and $1+\zeta^4=\lambda^{-1}\zeta^2$. 
Then
\begin{eqnarray*}
\mu_\zeta^k(\gamma)
=a_1(\kappa(x_1+\zeta^4))=\kappa(x_1+\zeta).
\end{eqnarray*}
We want to show the inequality $\stackrel{?}{<}$ in 
\begin{eqnarray}\label{7.7}
1&\stackrel{?}{<}&N_{\zeta:1}
\Bigl(\frac{\gamma}{b(\gamma)}\Bigr).
\end{eqnarray}
We calculate the right hand side as follows. The final step
uses that the Galois group 
$\textup{Gal}(\Q(\lambda)/\Q)=\{\id,\sigma\}$
contains the element $\sigma$ with $\sigma(\lambda)=-\lambda^{-1}=1-\lambda$.
Write $x_2:=\sigma(x_1)\in \Q(\lambda)-\{0\}$. Then 
\begin{eqnarray*}
N_{\zeta:1}\Bigl(\frac{\gamma}{b(\gamma)}\Bigr)
= N_{\zeta:1}\Bigl(\frac{\mu_\zeta^k(\gamma)}{b(\gamma)}\Bigr)
=N_{\zeta:1}\Bigl(\frac{x_1+\zeta}{x_1+\zeta^4}\Bigr)\stackrel{\textup{\eqref{2.1}}}{=} 
N_{\lambda:1}\Bigl(|\frac{x_1+\zeta}{x_1+\zeta^4}|^2\Bigr)\\
= N_{\lambda:1}\Bigl(\frac{x_1^2+\lambda x_1+1}{x_1^2-\lambda x_1+1}\Bigr)
=\frac{(x_1^2+\lambda x_1+1)(x_2^2-\lambda^{-1}x_2+1)}
{(x_1^2-\lambda x_1+1)(x_2^2+\lambda^{-1}x_2+1)}.
\end{eqnarray*}
We will show that the last quotient is $>1$ for any 
$x_1\in[1,\lambda]$ and any $x_2\in\R$. 
Remark that $x^2\pm \lambda^{\pm 1} x+1>0$ for any $x\in\R$
because $|\lambda|<2$ and $|\lambda^{-1}|<1$. 
Define the functions $f:\R\to\R_{>0}$ and 
$\www{f}:\R\to\R_{>0}$ by
\begin{eqnarray*}
f(x):=\frac{x^2+\lambda x+1}{x^2-\lambda x+1}
\quad\textup{and}\quad
\www{f}(x):=\frac{x^2-\lambda^{-1}x+1}{x^2+\lambda^{-1}x+1}.
\end{eqnarray*}
Their derivatives $f'$ and $\www{f}'$ are
\begin{eqnarray*}
f'(x)=\frac{-2\lambda(x^2-1)}{(x^2-\lambda x+1)^2}
\quad\textup{and}\quad
\www{f}'(x)=\frac{2\lambda^{-1}(x^2-1)}{(x^2+\lambda^{-1}x+1)^2}.
\end{eqnarray*}
Therefore
\begin{eqnarray*}
&&\min_{x\in[1,\lambda]}f(x)=f(\lambda)=2\lambda^2+1=2\lambda+3,\\
&&\min_{x\in\R}\www{f}(x)= \www{f}(1)=\frac{2-\lambda^{-1}}{2+\lambda^{-1}}
=7-4\lambda,\\
&&\Bigl(\min_{x\in[1,\lambda]}f(x)\Bigr)\Bigl(\min_{x\in\R}\www{f}(x)\Bigr)
=(2\lambda+3)(7-4\lambda)=1+6(2-\lambda)>1.
\end{eqnarray*}
\eqref{7.7} is proved. 

{\bf The case $n\geq 2$:} The element 
$b(\gamma)=a_1^{-n}(\mu_\zeta^k(\gamma))\in S[\zeta,\zeta^4)$ 
can be written as 
\begin{eqnarray*}
b(\gamma)=a_1^{-n}(\mu_\zeta^k(\gamma))= \kappa(x_1+\zeta^4)&&\textup{with suitable }
\kappa\in\Z[\lambda]_{>0}\\
&&\textup{and } x_1\in \kappa^{-1}\Z[\lambda]\cap (0,\lambda],
\end{eqnarray*}
Then
\begin{eqnarray*}
\mu_\zeta^k(\gamma)
=a_1^n(\kappa(x_1+\zeta^4))
=\kappa(x_1+n\lambda+\zeta^4).
\end{eqnarray*}
We want to show the inequality $\stackrel{?}{<}$ in 
\begin{eqnarray}\label{7.8}
1&\stackrel{?}{<}&N_{\zeta:1}\Bigl(\frac{\gamma}{b(\gamma)}\Bigr).
\end{eqnarray}
We calculate the right hand side as in the case $n=1$. 
Write $x_2:=\sigma(x_1)\in \Q(\lambda)-\{0\}$. Then 
\begin{eqnarray}
&& N_{\zeta:1}\Bigl(\frac{\gamma}{b(\gamma)}\Bigr)
= N_{\zeta:1}\Bigl(\frac{\mu_\zeta^k(\gamma)}{b(\gamma)}\Bigr)
=N_{\zeta:1}\Bigl(\frac{x_1+n\lambda+\zeta^4}{x_1+\zeta^4}\Bigr)
\nonumber \\
&&\stackrel{\textup{\eqref{2.1}}}{=} 
N_{\lambda:1}\Bigl(|\frac{x_1+n\lambda+\zeta^4}{x_1+\zeta^4}|^2\Bigr)
= N_{\lambda:1}\Bigl(\frac{(x_1+n\lambda)^2-(x_1+n\lambda)\lambda+1}{x_1^2-x_1\lambda+1}\Bigr) \nonumber \\
&&= N_{\lambda:1}\Bigl(\frac{x_1^2+(2n-1)\lambda x_1+n(n-1)\lambda^2+1}{x_1^2-\lambda x_1+1}\Bigr)\nonumber \\
&&=f(x_1)\cdot \frac{x_2^2-(2n-1)\lambda^{-1}x_2+n(n-1)\lambda^{-2}+1}
{x_2^2+\lambda^{-1}x_2+1},\label{7.9}
\end{eqnarray}
where the function $f:\R\to\R$ is defined by
\begin{eqnarray*}
f(x):=\frac{x^2+(2n-1)\lambda x+n(n-1)\lambda^2+1}{x^2-\lambda x+1}.
\end{eqnarray*}
We will show that the product in \eqref{7.9} is $>1$ for any 
$x_1\in(0,\lambda]$ and any $x_2\in\R$. 
The derivative $f'$ of the function $f$ is 
\begin{eqnarray*}
f'(x)=\frac{-2n\lambda\bigl((x^2+(n-1)\lambda x-(\frac{n-1}{2}\lambda^2+1)\bigr)}{(x^2-\lambda x+1)^2}.
\end{eqnarray*}
$f'$ has two real zeros $y_1$ and $y_2$ with $y_2<0$ and $y_1\in(0,\lambda)$,
because $f'(y)<0$ for $y\ll 0$, $f'(0)>0,f'(\lambda)<0$. 
Therefore 
\begin{eqnarray*}
\min_{x\in[0,\lambda]}f(x)&=&\min(f(0),f(\lambda))
=\min(n(n-1)\lambda^2+1,n(n+1)\lambda^2+1)\\
&=&n(n-1)\lambda^2+1.
\end{eqnarray*}
For the proof of \eqref{7.8} it is sufficient to prove for any $x\in\R$ 
the inequality $\stackrel{?}{<}$ in
\begin{eqnarray*}
0&\stackrel{?}{<}&\Bigl(n(n-1)\lambda^2+1\Bigr)
\Bigl(x^2-(2n-1)\lambda^{-1}x+n(n-1)\lambda^{-2}+1\Bigr)\\
&-&\Bigl(x^2+\lambda^{-1}x+1\Bigr).
\end{eqnarray*}
Setting $y:=x-(n-\frac{1}{2})\lambda^{-1}$, the right hand side transforms into
\begin{eqnarray*}
\Bigl(n(n-1)\lambda^2\Bigr)y^2+\Bigl(-2n\lambda^{-1}\Bigr)y
+\Bigl((n(n-1)\lambda^2)(1-\frac{1}{4}\lambda^{-2})-n^2\lambda^{-2}\Bigr)\\
=:a_1y^2+b_1y+c_1.
\end{eqnarray*}
This is positive for any $y\in\R$ if and only if its discriminant
$(\frac{b_1}{2a_1})^2-\frac{c_1}{a_1}$ is negative. In the following
estimates $n\geq 2$ and $\lambda=\frac{1+\sqrt{5}}{2}$ are used. 
\begin{eqnarray*}
(\frac{b_1}{2a_1})^2-\frac{c_1}{a_1}
&=& -1+\frac{1}{4}\lambda^{-2}+\frac{n}{n-1}\lambda^{-4}+\frac{1}{(n-1)^2}\lambda^{-6}\\
&\stackrel{n\geq 2}{\leq}& -1+\frac{1}{4}\lambda^{-2}+2\lambda^{-4}+\lambda^{-6}
=-\frac{1}{2}-\frac{1}{4}\lambda^{-5}-\lambda^{-7}<0.
\end{eqnarray*}
Here we used $(-\lambda^{-1})^k=F_{k+1}-F_k\lambda$, where $(F_n)_{n\geq 0}$ is the
Fibonacci sequence with $F_0=0$, $F_1=1$, $F_{n+2}=F_n+F_{n+1}$.
\eqref{7.8} is proved. 
\hfill$\Box$

\section{From triangular matrices to even and odd monodromy
groups and even and odd vanishing cycles}\label{s8}
\setcounter{equation}{0}
\setcounter{figure}{0}

We have some good reason to call
the subgroup $\langle a_1,a_2\rangle\subset G_q^\C$
{\it odd monodromy group} or {\it odd Coxeter group of rank 2}
and to call the elements of the set $\Delta^{(1)}_q$ 
{\it odd vanishing cycles}.
In this section we will explain this reason.
We start with a rather general construction.

\begin{definition}\label{t8.1}
Fix $n\in\N$ and a subring $R\subset\C$ with $1\in R$
(e.g. $R\in\{\Z,\Q,\R,\C\}$ or $R=\Z[2\cos\frac{\pi}{q}]$
for some $q\in\Z_{\geq 3}$). 

(a) Define the set
\begin{eqnarray}\label{8.1}
T^{uni}_n(R):=\{S\in M_{n\times n}(R)\,|\, S_{ij}=0
\textup{ for }i>j,S_{ii}=1\}
\end{eqnarray}
of upper triangular $n\times n$-matrices with diagonal entries
equal to 1 and entries in $R$. 

(b) Fix $S\in T^{uni}_n(R)$. 
Consider a free $R$-module $H_R$ of rank $n$ with a fixed
$R$-basis $\uuuu{e}=(e_1,...,e_n)$. Define a
symmetric $R$-bilinear form $I^{(0)}:H_R\times H_R\to R$
and a skew-symmetric $R$-bilinear form 
$I^{(1)}:H_R\times H_R\to R$ by
\begin{eqnarray}\label{8.2}
I^{(k)}(\uuuu{e}^t,\uuuu{e})=S+(-1)^kS^t 
\quad\textup{for }k\in\{0;1\}.
\end{eqnarray}
Define an $R$-linear automorphism $M:H_R\to H_R$ by 
\begin{eqnarray}\label{8.3}
M(\uuuu{e})=\uuuu{e}\cdot S^{-1}S^t.
\end{eqnarray}
Define for $a\in R$ with $I^{(0)}(a,a)=2$ an automorphism
$s_a^{(0)}:H_R\to H_R$, and define for arbitrary $a\in R$
an automorphism $s_a^{(1)}:H_R\to H_R$ by 
\begin{eqnarray}\label{8.4}
s_a^{(k)}(b)=b-I^{(k)}(a,b)a\quad\textup{for }b\in H_R
\textup{ and }k\in\{0;1\}.
\end{eqnarray}
Define two groups $\Gamma^{(0)}\subset\Aut(H_R)$ and
$\Gamma^{(1)}\subset\Aut(H_R)$ by 
\begin{eqnarray}\label{8.5}
\Gamma^{(k)}:=\langle s_{e_1}^{(k)},...,s_{e_n}^{(k)}\rangle
\quad\textup{for }k\in\{0;1\}.
\end{eqnarray}
Define two sets $\Delta^{(0)}\subset H_R$ and 
$\Delta^{(1)}\subset H_R$ by 
\begin{eqnarray}\label{8.6}
\Delta^{(k)}:=\Gamma^{(k)}\{\pm e_1,...,\pm e_n\}.
\end{eqnarray}
\end{definition}

Before presenting context where these data arise,
here are basic properties of them.
The parts (a) to (c) in the following lemma 
are elementary (e.g. \cite[Lemma 2.2]{HL24}).
Part (d) uses the triangular shape of $S$
(e.g. \cite[Theorem 2.7]{HL24}). 

\begin{lemma}\label{t8.2}
Fix $n,R,S$ and $H_R$ as in Definition \ref{t8.1}. 

(a) $M$ respects $I^{(0)}$ and $I^{(1)}$. 

(b) $s_a^{(0)}$ respects $I^{(0)}$ and is a reflection.
It is $\id$ on $\{b\in H_R\,|\, I^{(0)}(a,b)=0\}$
and $-\id$ on $\langle a\rangle$. Especially 
$(s_a^{(0)})^2=\id$.

(c) $s_a^{(1)}$ respects $I^{(1)}$ and is a transvection
if $a\notin\Rad I^{(1)}$ and $\id$ if $a\in \Rad I^{(1)}$.
If $a\notin\Rad I^{(1)}$ it is $\id$ on
$\{b\in H_R\,|\, I^{(1)}(a,b)=0\}$ and has a single
$2\times 2$ Jordan block (with eigenvalue 1).
Furthermore 
\begin{eqnarray}\label{8.7}
(s_a^{(1)})^{-1}(b)=b+I^{(1)}(a,b)a
\quad\textup{for }b\in H_R.
\end{eqnarray}

(d) 
\begin{eqnarray}\label{8.8}
s_{e_1}^{(k)}...s_{e_n}^{(k)}=(-1)^{k+1}M
\quad\textup{for }k\in\{0;1\}.
\end{eqnarray}
\end{lemma}

\begin{remarks}\label{t8.3}
(i) These data arise for $k=0$ and $k=1$ in the theory of
isolated hypersurface singularities. 
There $R=\Z$, $H_\Z$ is the Milnor lattice of a singularity,
$M$ is its {\it monodromy}, $I^{(0)}$ and $I^{(1)}$ are
its {\it even} and {\it odd intersection form},
$\Gamma^{(0)}$ and $\Gamma^{(1)}$ are its {\it even} and
{\it odd monodromy group}, and the elements of 
$\Delta^{(0)}$ and $\Delta^{(1)}$ are its {\it even} and
{\it odd vanishing cycles}.
For an introduction to this see e.g. \cite[10.1]{HL24}.

(ii) Parts of these data arise for $k=0$ in the theory of
Coxeter groups (see e.g. \cite[Chapter 5]{Hu90}).
Here $R=\R$, but $S\in T^{uni}_n(\R)$ 
is chosen quite special, with 
\begin{eqnarray*}
S_{ij}&=&-2\cos\frac{\pi}{q_{ij}}\quad\textup{for some }
q_{ij}\in\Z_{\geq 2}\\
\textup{or }S_{ij}&\in& \R_{\leq -2}.
\end{eqnarray*}
Define $I:=\{(i,j)\,|\, 1\leq i<j\leq n, 
S_{ij}=-2\cos\frac{\pi}{q_{ij}}\}.$
Then $\Gamma^{(0)}$ is a realization of a {\it Coxeter group}.
A Coxeter group is a group with a presentation
\begin{eqnarray}\label{8.9}
\langle \www{r}_1,....,\www{r}_n\,|\, 
\www{r}_1^2=...=\www{r}_n^2=e, 
(\www{r}_i\www{r}_j)^{q_{ij}}=e\textup{ for }(i,j)\in I\}.
\end{eqnarray}
The isomorphism from $\Gamma^{(0)}$ to this group 
is given by $s_{e_i}^{(0)}\to \www{r}_i$.
This fact is classical (see e.g. \cite[5.3 and 5.4]{Hu90})
if $S_{ij}=-2$ for $(i,j)\notin I$. It was generalized
to the case $S_{ij}\in \R_{\leq -2}$ for $(i,j)\notin I$ by
Vinberg \cite[Proposition 6, Theorem 1, Theorem 2,
Proposition 17]{Vi71}.

(iii) Comparing (i) and (ii), it is natural to ask whether
the theory of Coxeter groups has an odd variant, i.e. a variant
for $k=1$. 
In the theory of Coxeter groups the building blocks
are the rank 2 Coxeter groups. Therefore it is natural to study
for a matrix 
$$S=\begin{pmatrix}1&-\lambda\\0&1\end{pmatrix} 
\quad\textup{with }\lambda=2\cos\frac{\pi}{q}
\textup{ and }q\in\Z_{\geq 3}$$
the odd monodromy group $\Gamma^{(1)}$. 
This leads us back to the main subject of this paper. 
The next theorem puts things together.
\end{remarks}

\begin{theorem}\label{t8.4}
Fix the data from the introduction, a natural number $q\geq 3$,
$\zeta=e^{2\pi i/(2q)}$ and $\lambda=\zeta+\oooo{\zeta}
=2\cos\frac{\pi}{q}\in [1,2)$. Consider the rings 
$R:=\Z[\lambda]$ and $\Z[\zeta]\stackrel{2:1}{\supset}\Z[\lambda]$ 
and the matrix
\begin{eqnarray}\label{8.10}
S=\begin{pmatrix}1&-\lambda\\0&1\end{pmatrix}\in T^{uni}_2(\Z[\lambda]).
\end{eqnarray}
We choose as free $\Z[\lambda]$-module of rank 2 
\begin{eqnarray}\label{8.11}
H_R:=\Z[\zeta]
\quad\textup{with }\Z[\lambda]\textup{-basis}\quad 
\uuuu{e}=(e_1,e_2)=(1,\zeta^{q-1}).
\end{eqnarray}
Then $H_\R=\C$. 

(a) The matrices of $I^{(0)}$ and $I^{(1)}$ with respect to
the $\R$-bases $\uuuu{e}$ and $(1,i)$ of $H_\R=\C$ are 
\begin{eqnarray*}
I^{(0)}(\uuuu{e}^t,\uuuu{e})=S+S^t
=\begin{pmatrix}2&-\lambda\\-\lambda&2\end{pmatrix}, \quad
I^{(1)}(\uuuu{e}^t,\uuuu{e})=S-S^t
=\begin{pmatrix}0&-\lambda\\\lambda&0\end{pmatrix},\\
I^{(0)}(\begin{pmatrix}1\\i\end{pmatrix},(1,i))
=\begin{pmatrix}2&0\\0&2\end{pmatrix},\quad
I^{(1)}(\begin{pmatrix}1\\i\end{pmatrix},(1,i))
=\begin{pmatrix}0&-\frac{\lambda}{\sqrt{1-\frac{\lambda^2}{4}}}
\\ \frac{\lambda}{\sqrt{1-\frac{\lambda^2}{4}}}&0\end{pmatrix}.
\end{eqnarray*}
So, $I^{(0)}$ is two times the standard scalar product
on $\C=\R\cdot 1\oplus\R\cdot i$.

(b) Recall from \eqref{5.6} and \eqref{5.7}
$\Psi(Q)=\mu_\zeta$, $\Psi(V)=v$, $\Psi(A_1)=a_1$, 
$\Psi(A_2)=a_2$ with  
\begin{eqnarray}\label{8.12}
v^2=-\id=\mu_{-1}\quad\textup{and}\quad 
a_1v=\mu_\zeta=va_2.
\end{eqnarray}
We have 
\begin{eqnarray}\label{8.13}
a_1=s_{e_1}^{(1)},\quad a_2=s_{e_2}^{(1)},\\
\Gamma^{(1)}=\langle a_1,a_2\rangle \subset G_q^\C.\label{8.14}
\end{eqnarray}

(c) Define\begin{eqnarray}
\www{q}:=
\left\{\begin{array}{ll}q&\textup{ if }q\equiv 0(4),\\
q/2&\textup{ falls }q\equiv 2(4),\\
2q&\textup{ falls }q\equiv 1(2),\end{array}\right. 
\label{8.15}
\end{eqnarray}
Then
\begin{eqnarray}\label{8.16}
a_1a_2&=&\mu_{-\zeta^2}\quad\textup{and}\quad 
(a_1a_2)^{\www{q}}=\id.
\end{eqnarray}
If $q$ is even, then
\begin{eqnarray}\label{8.17}
(a_1a_2)^{q/2}=(a_2a_1)^{q/2}=\left\{\begin{array}{ll}
\mu_{-1}&\textup{if }q\equiv 0(4),\\
\id&\textup{if }q\equiv 2(4).\end{array}\right. 
\end{eqnarray}
If $q$ is odd, then
\begin{eqnarray}\label{8.18}
(a_1a_2a_1...)_{[q\textup{ factors}]}
&=&(a_2a_1a_2...)_{[q\textup{ factors}]}
=(-1)^{(q-1)/2}v,\\
(a_1a_2)^{(q+1)/2}&=& (-1)^{(q-1)/2}\mu_\zeta.\label{8.19}
\end{eqnarray}

(d) If $q$ is odd then 
\begin{eqnarray}\label{8.20}
\Gamma^{(1)}=G_q^\C.
\end{eqnarray}
\end{theorem}

{\bf Proof:}
(a) The first line of part (a) is the definition 
\eqref{8.2} of $I^{(0)}$ and $I^{(1)}$. The second line
follows with the base change matrix 
\begin{eqnarray*}
B:=\begin{pmatrix}1&-\lambda/2\\0&\sqrt{1-\lambda^2/4}\end{pmatrix}
\quad\textup{with}\quad 
\uuuu{e}=(1,\zeta^{q-1})=(1,i)B\\
\textup{and}\quad 
I^{(k)}(\uuuu{e}^t,\uuuu{e})
=B^t\cdot I^{(k)}(1,i)^t,(1,i))\cdot B.
\end{eqnarray*}

(b) \eqref{8.12} is \eqref{5.9}. \eqref{8.13}
follows if one writes out the matrix of $s_{e_i}^{(1)}$ 
with respect to $\uuuu{e}=(1,\zeta^{q-1})$, which is $A_i$.

(c) $a_1v=\mu_\zeta=va_2$ and $v^2$ give
$a_1a_2=-\mu_\zeta^2=\mu_{-\zeta^2}$. 
It implies $(a_1a_2)^{\www{q}}=\id$, \eqref{8.17}
and \eqref{8.19}. Together \eqref{8.19} and \eqref{8.12}
imply \eqref{8.18}. 

(d) Because of \eqref{8.19} $\Gamma^{(1)}$ contains
$\pm \mu_\zeta$, so $\pm v$, so $v$ and $\mu_\zeta$. 
\hfill$\Box$

\bigskip
For odd $q$ we know $\Gamma^{(1)}=G_q^\C$ rather well. 
For even $q$ we first consider the quotient group
$\oline{\Gamma^{(1),mat}}\subset G_q\subset PSL_2(\Z[\lambda])$
of the matrix version 
$\Gamma^{(1),mat}:=\langle A_1,A_2\rangle$ of $\Gamma^{(1)}$.
Recall the Remarks \ref{t4.1} on triangle groups and the
points $p_1=\infty$, $p_2=\zeta$, $p_3=i$ and $p_4=0$
in Theorem \ref{t4.2}.

\begin{lemma}\label{t8.5}
Suppose that $q$ is even.

(a) The group $\oline{\Gamma^{(1),mat}}=\langle \oline{A_1},
\oline{A_2}\rangle$ is the triangle group for the hyperbolic
triangle with vertices $p_1,p_2,p_4$. 
It is a normal subgroup of index 2 in the triangle group
$G_q$. The map $\oline{A_1}\to\www{a_1}$, 
$\oline{A_2}\to\www{a_2}$ extends to an isomorphism from
$\oline{\Gamma^{(1),mat}}$ to the group with presentation
\begin{eqnarray}\label{8.21}
\langle \www{a_1},\www{a_2}\,|\, (\www{a_1}\www{a_2})^{q/2}=e
\rangle,
\end{eqnarray}
so $\oline{\Gamma^{(1),mat}}$ is isomorphic to the free product
$\Z\star \Z_{q/2}$. 

(b) A fundamental domain $\FF_2\subset\H$ for the action of
$\oline{\Gamma^{(1),mat}}$ is the set $\FF_2$ with 
\begin{eqnarray}\nonumber
\oooo{\FF_2}&:=&
\textup{the degenerate hyperbolic quadrangle with vertices}\\
&& \hspace*{1cm}\infty,\zeta,0,\zeta^{q-1},\label{8.22}
\end{eqnarray}
and boundary as prescribed by \eqref{4.7} and \eqref{4.8}
(with $\FF$ replaced by $\FF_2$). 
The stabilizers of the
points $\zeta$, $\infty$ and $0$ are 
\begin{eqnarray}\label{8.23}
\Stab_{{\oline{\Gamma^{(1),mat}}}}(\zeta)
=\langle{\oline{A_1A_2}}\rangle,\\
\Stab_{{\oline{\Gamma^{(1),mat}}}}(\infty)
=\langle{\oline{A_1}}\rangle,\ 
\Stab_{{\oline{\Gamma^{(1),mat}}}}(0)
=\langle{\oline{A_2}}\rangle.\label{8.24}
\end{eqnarray}
For any $z\in\H$, the orbit ${\oline{\Gamma^{(1),mat}}}\{z\}$
intersects $\FF_2$ in one point. 
If ${\oline{C}}\in{\oline{\Gamma^{(1),mat}}}$
and $z\in\FF_2$ with ${\oline{C}}(z)\in\FF_2$ then 
${\oline{C}}(z)=z$ and 
\begin{eqnarray}\label{8.25}
\Bigl({\oline{C}}=\id\Bigr)
\textup{ or }\Bigl(z=\zeta,\ {\oline{C}}
\in \langle{\oline{A_1A_2}}\rangle\Bigr).
\end{eqnarray}
\end{lemma}

{\bf Proof:}
(a) The proof of Theorem \ref{t4.2} showed how
$\oline{A_1}$ and $\oline{A_2}$ act on $\H$. 
They are parabolic with fixed point $\infty$ respectively
$0$. Their product $\oline{A_1A_2}=\oline{Q^2}$ is elliptic
with fixed point $\zeta$ and rotation angle $-\frac{4\pi}{q}$.
As in the proof of Theorem \ref{t4.2} one sees that
$\oline{\Gamma^{(1),mat}}= 
\langle \oline{A_1},\oline{A_2}\rangle$
is the claimed triangle group. Together with the Remarks 
\ref{t4.1}, this shows also the presentation \eqref{8.21}
of the group $\oline{\Gamma^{(1),mat}}$. 

(b) This follows exactly as Corollary \ref{t4.3}.
\hfill$\Box$

\section{The odd Coxeter like groups of rank 2}\label{s9}
\setcounter{equation}{0}
\setcounter{figure}{0}

This section continues the study of the odd monodromy
group $\Gamma^{(1)}=\langle a_1,a_2\rangle$ in the situation
of Theorem \ref{t8.4}. This is the odd analog of the
Coxeter group of rank 2, which is of type $I_2(q)$. 
Especially, we give presentations of $G_q^\C$ 
(Lemma \ref{t9.1} (a)) and $\Gamma^{(1)}\subset G_q^\C$
(Theorem \ref{t9.2} (b)). Theorem \ref{t9.2}
(a) will show that the set $\Delta^{(1)}_q=G_q^\C\{1\}$
from the introduction coincides here with the
set $\Delta^{(1)}$ from Definition \ref{t8.1}. 
It treats the subgroup
$\Gamma^{(1)}\subset G_q^\C$ in all cases.
We have the cases $q\equiv 1(2)$, $q\equiv 0(4)$
and $q\equiv 2(4)$.
At the end of this section we compare $\Gamma^{(1)}$
with the corresponding rank 2 Coxeter group.

The next lemma uses Corollary \ref{t4.3} to give a 
presentation of $G_q^\C$ for all $q$, and it uses
Lemma \ref{t8.5} (b) to give a presentation of
$\Gamma^{(1)}$ in the case $q\equiv 2(4)$.

\begin{lemma}\label{t9.1}
(a) Consider the group $G_q^\C=\langle v,\mu_\zeta\rangle
=\langle v,a_1\rangle=\langle a_1,\mu_\zeta\rangle$. 
The map $v\mapsto \whh{v}$, $\mu_\zeta\mapsto\whh{q}$,
extends to an isomorphism from 
$G_q^\C$ to the group with presentation
\begin{eqnarray}\label{9.1}
\langle \whh{v},\whh{q}\,|\, 
\whh{v}^2=\whh{q}^q,\whh{v}^4=e\rangle.
\end{eqnarray}

(b) Suppose that $q\equiv 2(4)$. Consider the group 
$\Gamma^{(1)}=\langle a_1,a_2\rangle$. The map 
$a_1\mapsto\www{a_1}$, $a_2\mapsto\www{a_2}$, 
extends to an isomorphism from 
$\Gamma^{(1)}$ to the group with presentation
\begin{eqnarray}\label{9.2}
\langle \www{a_1},\www{a_2}\,|\, 
(\www{a_1}\www{a_2})^{q/2}=e\rangle,
\end{eqnarray}
so $\langle a_1,a_2\rangle \cong \Gamma^{(1)} 
\cong {\oline{\Gamma^{(1),mat}}}\cong\Z\star\Z_{q/2}$,
and $\mu_{-1}\notin\langle a_1,a_2\rangle$.
\end{lemma}

{\bf Proof:}
(a) Because of $v^2=\mu_{-1}=\mu_\zeta^q$,
the map $\whh{v}\mapsto v,\whh{q}\mapsto\mu_\zeta$
extends to a homomorphism from the group 
$\langle \whh{v},\whh{q}\rangle$
with presentation in \eqref{9.1} to the group
$\langle v,\mu_\zeta\rangle$. 
The map $\whh{v}\to\www{v},\whh{q}\mapsto\www{q}$
extends to a surjective group homomorphism from the group
in \eqref{9.1} to the group in \eqref{4.5}, whose kernel
$\{\whh{v}^2,e\}$ has order 2. 
On the other hand, the map 
$V\mapsto \www{v},Q\mapsto \www{q}$, extends because
of Theorem \ref{t4.2} (c) also to a surjective group homomorphism
with kernel of order 2. A diagram for this situation:
\begin{eqnarray*}
\langle \whh{v},\whh{q}\rangle 
\to
\langle v,\mu_\zeta\rangle 
\cong
\langle V,Q\rangle
\stackrel{2:1}{\longrightarrow}
\langle {\oline{V}},{\oline{Q}}\rangle
\cong 
\langle \www{v},\www{q}\rangle
\stackrel{2.1}{\longleftarrow}
\langle \whh{v},\whh{q}\rangle.
\end{eqnarray*}
Here the first and last group are equal.
Therefore the homomorphism $\langle \whh{v},\whh{q}\rangle
\to \langle v,\mu_\zeta\rangle$ is an isomorphism.

(b) Because of $(a_1a_2)^{q/2}=\id$ in \eqref{8.16},
the map $\www{a_1}\mapsto a_1,\www{a_2}\mapsto a_2$
extends to a homomorphism from the group 
$\langle \www{a_1},\www{a_2}\rangle$
with presentation in \eqref{9.2} to the group 
$\langle a_1,a_2\rangle$.  
The groups in \eqref{9.2} and in \eqref{8.21}
coincide and are because of Lemma \ref{t8.5} (a) 
isomorphic to the group
$\langle {\oline{A_1}},{\oline{A_2}}\rangle$. 
Of course, we have the canonical homomorphism
$\langle A_1,A_2\rangle  \to 
\langle {\oline{A_1}},{\oline{A_2}}\rangle$. 
A diagram for this situation:
\begin{eqnarray*}
\langle \www{a_1},\www{a_2}\rangle
\to
\langle a_1,a_2\rangle
\cong
\langle A_1,A_2\rangle
\to 
\langle {\oline{A_1}},{\oline{A_2}}\rangle
\cong
\langle \www{a_1},\www{a_2}\rangle.
\end{eqnarray*}
Here the first and last group are equal.
Therefore all homomorphisms in this diagram are
isomorphisms. 
Especially $-E_2\notin \langle A_1,A_2\rangle$
and $\mu_{-1}\notin\langle a_1,a_2\rangle$. 
\hfill$\Box$ 

\bigskip
In Theorem \ref{t9.2} we work with the set
$\Delta^{(1)}_q:=G_q^\C\{1\}$ from the introduction.
Part (a) of Theorem \ref{t9.2} 
shows that it coincides with the set
$\Delta^{(1)}$ in Definition \ref{t8.1}.

\begin{theorem}\label{t9.2}
(a) (i) Suppose that $q\equiv 1(2)$. Then
\begin{eqnarray*}
\Gamma^{(1)} = G_q^\C\quad
\textup{and}\quad \mu_{-1}\in \Gamma^{(1)}.
\end{eqnarray*}

(ii) Suppose that $q\equiv 0(4)$. Then
$\Gamma^{(1)}$ is a normal subgroup of index 2 in $G_q^\C$, and 
\begin{eqnarray*}
\mu_{-1}\in \Gamma^{(1)},\quad
\langle \mu_\zeta\rangle\cap\Gamma^{(1)}
=\langle \mu_{\zeta^2}\rangle\textup{ has index 2 in }
\langle \mu_\zeta\rangle,\\
\Delta^{(1)}_q=\Gamma^{(1)}\{1\}\ \dot\cup\  
\Gamma^{(1)}\{\zeta^{q-1}\}.
\end{eqnarray*}

(iii) Suppose that $q\equiv 2(4)$. Then
$\Gamma^{(1)}$ is a normal subgroup of index 4 in $G_q^\C$, and 
\begin{eqnarray*}
\mu_{-1}\notin \Gamma^{(1)},\quad
\langle \mu_\zeta\rangle\cap\Gamma^{(1)}
=\langle \mu_{\zeta^4}\rangle\textup{ has index 4 in }
\langle \mu_\zeta\rangle,\\
\Delta^{(1)}_q=\Gamma^{(1)}\{1\}\ \dot\cup\  
\Gamma^{(1)}\{-1\}\ \dot\cup\  
\Gamma^{(1)}\{\zeta^{q-1}\}\ \dot\cup\  
\Gamma^{(1)}\{-\zeta^{q-1}\}.
\end{eqnarray*}

(b) Now consider any $q\in\Z_{\geq 2}$. 
The map $a_1\mapsto \whh{a_1}$, $a_2\mapsto \whh{a_2}$
extends to an isomorphism from the group 
$\Gamma^{(1)}$ to the group with presentation
\begin{eqnarray}\label{9.3}
\langle \whh{a_1},\whh{a_2}\,|\, 
(\whh{a_1}\whh{a_2}\whh{a_1}...)_{[q\textup{ factors}]}
=(\whh{a_2}\whh{a_1}\whh{a_2}...)_{[q\textup{ factors}]},
(\whh{a_1}\whh{a_2})^{\www{q}}=e\rangle
\end{eqnarray}
with $\www{q}$ as in \eqref{8.15}. 
\end{theorem}

{\bf Proof:}
(a) (i) See Theorem \ref{t8.4} (d). 

(ii) The part $(a_1a_2)^{q/2}=\mu_{-1}$ of \eqref{8.17}
shows $\mu_{-1}\in\langle a_1,a_2\rangle$ and
$-E_2\in\langle A_1,A_2\rangle$. 
Therefore $\Gamma^{(1),mat}$ and $G_q^{mat}$
are the full preimages in $SL_2(\R)$ of the subgroups
$\oline{\Gamma^{(1),mat}}$ and $G_q$ of $PSL_2(\R)$.
Therefore the index of $\Gamma^{(1),mat}$ in 
$G_q^{mat}$ is the same as the index of 
$\oline{\Gamma^{(1),mat}}$ in $G_q$, which is 2
by Lemma \ref{t8.5} (a). 

$1\in\C$ has the same stabilizer $\langle a_1\rangle$
in $G_q^\C$ and in its index 2 subgroup $\Gamma^{(1)}$.
Therefore $\Delta^{(1)}_q=G_q^\C\{1\}$ splits into two 
$\Gamma^{(1)}$-orbits. Because of 
$$G_q^\C=\Gamma^{(1)}\ \dot\cup\  \mu_{\zeta^{q-1}}\Gamma^{(1)}$$
by $\langle\mu_\zeta\rangle\cap\Gamma^{(1)}
=\langle \mu_{\zeta^2}\rangle$, 
these are the $\Gamma^{(1)}$-orbits of $1$ and of $\zeta^{q-1}$.

(iii) $G_q^{mat}$ is the full preimage in
$SL_2(\R)$ of $G_q$,
but $\Gamma^{(1),mat}$ is because of 
$-{\bf 1}_2\notin \langle A_1,A_2\rangle$ (Lemma \ref{t9.1} (b))
an isomorphic preimage in $SL_2(\R)$ of 
$\oline{\Gamma^{(1),mat}}$.
Therefore the index of $\Gamma^{(1),mat}$
in $G_q^{mat}$ is 4. It is a normal subgroup because
$\Gamma^{(1)}$ is a normal subgroup of 
$G_q^\C$, and this holds because 
$\langle a_1,a_2,v\rangle= G_q^\C$ and 
\begin{eqnarray*}
v^{-1} a_1v&=& a_2, \quad v^{-1}a_2v=va_2v^{-1}=a_1.
\end{eqnarray*}

$1\in\C$ has the same stabilizer $\langle a_1\rangle$
in $G_q^\C$ and in its index 2 subgroup $\Gamma^{(1)}$.
Therefore $\Delta^{(1)}_q=G_q^\C\{1\}$ splits into four 
$\Gamma^{(1)}$-orbits. Because of 
$$G_q^\C=\Gamma^{(1)}\ \dot\cup\  \mu_{-1}\Gamma^{(1)}
\ \dot\cup\  \mu_{\zeta^{q-1}}\Gamma^{(1)}\ \dot\cup\  
\mu_{-\zeta^{q-1}}\Gamma^{(1)}$$
by $\langle\mu_\zeta\rangle\cap\Gamma^{(1)}
=\langle \mu_{\zeta^4}\rangle$, 
these are the $\Gamma^{(1)}$-orbits of $1$, $-1$, $\zeta^{q-1}$
and $-\zeta^{q-1}$.

(b) {\bf The case $q\equiv 2(4)$:} 
Then $\www{q}=\frac{q}{2}$, and the relation
$(\whh{a_1}\whh{a_2})^{\www{q}}=e$ shows that both sides
in the relation 
$(\whh{a_1}\whh{a_2}\whh{a_1}...)_{[q\textup{ factors}]}
=(\whh{a_2}\whh{a_1}\whh{a_2}...)_{[q\textup{ factors}]}$
are equal to $e$, so this relation holds if
$(\whh{a_1}\whh{a_2})^{\www{q}}=e$ holds.
Therefore Lemma \ref{t9.1} (b) implies part (b) in this case.

{\bf The case $q\equiv 0(4)$:} In this case $\www{q}=q$, and 
one can write \eqref{9.3} also as
\begin{eqnarray}\label{9.4}
\langle \whh{a_1},\whh{a_2}\,|\, 
(\whh{a_1}\whh{a_2})^{q/2}
=(\whh{a_2}\whh{a_1})^{q/2}, 
(\whh{a_1}\whh{a_2})^{q}=e\rangle.
\end{eqnarray}
Because of $(a_1a_2)^{q/2}=(a_2a_1)^{q/2}=\mu_{-1}$ in 
\eqref{8.17}, the map 
$\whh{a_1}\mapsto a_1$, $\whh{a_2}\mapsto a_2$, extends
to a homomorphism from the group 
$\langle\whh{a_1},\whh{a_2}\rangle$ with presentation in 
\eqref{9.4} to the group $\langle a_1,a_2\rangle$. 
The map 
$\whh{a_1}\mapsto\www{a_1},\whh{a_2}\mapsto\www{a_2}$
extends to  surjective group homomorphism from the group
in \eqref{9.4} to the group in \eqref{9.2},
whose kernel $\{(\whh{a_1}\whh{a_2})^{q/2},e\}$ has order 2. 
On the other hand, $\langle\oline{A_1},\oline{A_2}\rangle
\cong\langle\www{a_1},\www{a_2}\rangle$ by 
Lemma \ref{t9.1} (b), and
the map $A_1\mapsto\oline{A_1}\mapsto\www{a_1}$,
$A_2\mapsto\oline{A_2}\mapsto\www{a_2}$ is because of 
$-E_2\in\langle A_1,A_2\rangle$ surjective with kernel
of order 2. A diagram for this situation:
\begin{eqnarray*}
\langle \whh{a_1},\whh{a_2}\rangle
\to
\langle a_1,a_2\rangle
\cong
\langle A_1,A_2\rangle
\stackrel{2:1}{\longrightarrow}
\langle {\oline{A_1}},{\oline{A_2}}\rangle
\cong
\langle\www{a_1},\www{a_2}\rangle
\stackrel{2:1}{\longleftarrow}
\langle \whh{a_1},\whh{a_2}\rangle.
\end{eqnarray*}
Here the first and last group are equal.
Therefore the map 
$\langle \whh{a_1},\whh{a_2}\rangle\to
\langle a_1,a_2\rangle$ is an isomorphism.

{\bf The case $q\equiv 1(2)$:} Then $\www{q}=2q$. 
Because of part (a)(i) and Lemma \ref{t9.1} (a), 
it is sufficient to show the following three points:
(1) The map
\begin{eqnarray*}
&&\whh{v}\mapsto (\whh{a_1}\whh{a_2})^{q(q-1)/2}
(\whh{a_1}\whh{a_2}\whh{a_1}...)_{[q \textup{ factors}]},\\
&&\whh{q}\mapsto 
(\whh{a_1}\whh{a_2})^{q(q-1)/2}(\whh{a_1}\whh{a_2})^{(q+1)/2},
\end{eqnarray*}
extends to a homomorphism from the group in \eqref{9.1}
to the group in \eqref{9.3}. 
(2) The map 
\begin{eqnarray*}
\whh{a_1}\mapsto \whh{q}\whh{v}^{-1},&&
\whh{a_2}\mapsto \whh{v}^{-1}\whh{q}
\end{eqnarray*}
extends to a homomorphism from the group in \eqref{9.3}
to the group in \eqref{9.1}. 
(3) The composition of the map in (2) with the map in (1) 
gives the identity on the group in \eqref{9.3}.

All three points (1), (2) and (3) are easy to see.
Here observe that in \eqref{9.1} $\whh{v}^2=\whh{q}^q$
commutes with $\whh{v}$ and $\whh{q}$ and that in \eqref{9.3}
$(\whh{a_1}\whh{a_2})^q$ commutes with $\whh{a_1}$ and $\whh{a_2}$.
\hfill$\Box$ 

\bigskip
The cases $q=3$, $q=4$ and $q=6$ are three cases
with $q\equiv 1(2)$, $q\equiv 0(4)$ and $q\equiv 2(4)$.
The following examples add to Theorem \ref{t6.1} the
group $\Gamma^{(1),mat}$ and the way how $\Delta^{(1)}_q$ splits
into $\Gamma^{(1)}$-orbits in the cases $q=4$ and $q=6$.

\begin{examples}\label{t9.3}
(i) Let $q=4$. Then $\lambda=\sqrt{2}$ and 
\begin{eqnarray*}
\Gamma^{(1),mat}&=&
\{\begin{pmatrix}a&b\sqrt{2}\\c\sqrt{2}&d\end{pmatrix}
\,|\, a,b,c,d\in\Z,ad-2bc=1,\\
&&\hspace*{3cm}a\equiv 1(2),d\equiv 1(2)\},\\
\Gamma^{(1)}\{1\}&=&
\{a+c\sqrt{2}\zeta^3\,|\, a,c\in\Z,\gcd(a,2c)=1\},\\
\Gamma^{(1)}\{\zeta^3\}&=&
\{b\sqrt{2}+d\zeta^3\,|\, b,d\in\Z,\gcd(2b,d)=1\},
\nonumber\\
\Delta^{(1)}_q&=&\Gamma^{(1)}\{1\}\ \dot\cup\ 
\Gamma^{(1)}\{\zeta^3\}.
\end{eqnarray*}

(c) Let $q=6$. Then $\lambda=\sqrt{3}$ and 
\begin{eqnarray*}
\Gamma^{(1),mat}&=&
\{\begin{pmatrix}a&b\sqrt{3}\\c\sqrt{3}&d\end{pmatrix}
\,|\, a,b,c,d\in\Z,ad-3bc=1,\\
&&\hspace*{3cm}a\equiv 1(3),d\equiv 1(3)\},\\
\Gamma^{(1)}\{1\}&=&
\{a+c\sqrt{3}\zeta^5\,|\, a,c\in\Z,\gcd(a,c)=1,a\equiv 1(3)\},\\
\Gamma^{(1)}\{-1\}&=&
\{a+c\sqrt{3}\zeta^5\,|\, a,c\in\Z,\gcd(a,c)=1,a\equiv 2(3)\},\\
\Gamma^{(1)}\{\zeta^5\}&=&
\{b\sqrt{3}+d\zeta^5\,|\, b,d\in\Z,\gcd(b,d)=1,d\equiv 1(3)\},\\
\Gamma^{(1)}\{-\zeta^5\}&=&
\{b\sqrt{3}+d\zeta^5\,|\, b,d\in\Z,\gcd(b,d)=1,d\equiv 2(3)\},\\
\Delta^{(1)}_q&=&\Gamma^{(1)}\{1\}\ \dot\cup\  
\Gamma^{(1)}\{-1\}\ \dot\cup\  
\Gamma^{(1)}\{\zeta^5\}\ \dot\cup\ 
\Gamma^{(1)}\{-\zeta^5\}.
\end{eqnarray*}
\end{examples}

Definition \ref{t8.1} gives in the setting of Theorem \ref{t8.4}
two monodromy groups, the even monodromy group $\Gamma^{(0)}$
and the odd monodromy group $\Gamma^{(1)}$. 
Our claim that $\Gamma^{(1)}$ is the odd analog of a rank 2
Coxeter group builds on the fact that $\Gamma^{(0)}$
is a rank 2 Coxeter group. The following Remarks
discuss this.

\begin{remarks}\label{t9.4}
Consider the situation of Theorem \ref{t8.4}.

(i) Write
\begin{eqnarray*}
r_1:=s_{e_1}^{(0)},\quad r_2:=s_{e_2}^{(0)},\\
\textup{so that }\Gamma^{(0)}=\langle r_1,r_2\rangle.
\end{eqnarray*}
Recall that $I^{(0)}$ equips $\C=\R\cdot 1\oplus\R\cdot i$
with two times the standard scalar product.
$r_1$ and $r_2$ are by definition the reflections along the lines
orthogonal to $e_1=1$ and $e_2=\zeta^{q-1}$. Their product is
$$r_1r_2=\mu_{\zeta^2},$$
so the rotation by the angle $\frac{2\pi}{q}$. 
Let $\sigma_c$ for $c\in S^1$ denote the reflection along
the line orthogonal to $c$, so that $r_1=\sigma_1$
and $r_2=\sigma_{\zeta^{q-1}}$. Then
\begin{eqnarray*}
\Gamma^{(0)}=\{\sigma_{\zeta^k},\mu_{\zeta^{2k}}\,|\, 
k\in\{0,1,...,q-1\}\}
\end{eqnarray*}
so it is indeed the rank 2 Coxeter group of type $I_2(q)$
and isomorphic to the dihedral group $D_{2q}$.

(ii) It is well known that the map $r_1\mapsto \www{r_1}$,
$r_2\mapsto\www{r_2}$, extends to an isomorphism from
$\Gamma^{(0)}$ to the group with presentation
\begin{eqnarray}\label{9.5}
\langle \www{r_1},\www{r_2}\,|\, \www{r_1}^2=\www{r_2}^2=e,\
(\www{r_1}\www{r_2})^q=e\rangle.
\end{eqnarray}

(iii) From (ii) one obtains immediately the following 
closely related presentation of $\Gamma^{(0)}$.
The map $r_1\mapsto \whh{r_1}$,
$r_2\mapsto\whh{r_2}$, extends to an isomorphism from
$\Gamma^{(0)}$ to the group with presentation
\begin{eqnarray}\label{9.6}
\langle \whh{r_1},\whh{r_2}\,|\, \whh{r_1}^2=\whh{r_2}^2=e,\
(\whh{r_1}\whh{r_2}\whh{r_1}...)_{[q\textup{ factors}]}
=(\whh{r_2}\whh{r_1}\whh{r_2}...)_{[q\textup{ factors}]}\rangle.
\end{eqnarray}
This presentation has some similarity with the presentation
of $\Gamma^{(1)}$ in Theorem \ref{t9.2} (b). 
But it is not clear how this similarity could generalize
to higher rank. On the contrary, see the next Remarks \ref{t9.5}.

(iv) The set $\Delta^{(0)}$ of even vanishing cycles in the
setting of Theorem \ref{t8.4} is simply the set of
$2q$-th unit roots 
$\UR_{2q}=\{\zeta^k\,|\, k\in\{0,1,...,2q-1\}\}$,
\begin{eqnarray}\label{9.7}
\Delta^{(0)}=\UR_{2q}.
\end{eqnarray}
\end{remarks}

\begin{remarks}\label{t9.5}
One motivation to consider the groups 
$\Gamma^{(1)}\subset G_q^\C$ in Theorem \ref{t8.4}
is their analogy to the rank 2 Coxeter groups, which are the 
building blocks of arbitrary Coxeter groups.
If one wants to create a theory of odd analoga of
Coxeter groups, one has to start with the rank 2 cases.

Unfortunately, already in rank 3 optimistic expectations are
not met. The matrix
$$S=\begin{pmatrix}1&-2&-2\\0&1&-2\\0&0&1\end{pmatrix}
\in T^{uni}_3(\Z)$$
gives rise to the free Coxeter group $\Gamma^{(0)}$ with
three generators by \eqref{8.9}.
One might hope that $\Gamma^{(1)}$ from this matrix would be the
free group with three generators. But that is not the case.
By Theorem 6.18 (e) (ii) in \cite{HL24},
the group $\Gamma^{(1)}$ has a normal subgroup isomorphic to
$(\Z^2,+)$, and the quotient is isomorphic to the
product of $\{\pm 1\}$ with the free group with two
generators.

For quite many other matrices $S\in T^{uni}_3(\Z)$ the
group $\Gamma^{(1)}$ is by \cite[Theorem 6.18 (f)]{HL24} 
indeed the free group with three generators, but not for 
$S$ above. This leaves us puzzled. 
\end{remarks}

\section{More on the literature}\label{s10}
\setcounter{equation}{0}
\setcounter{figure}{0}

\subsection{$\lambda$-continued fractions}\label{s10.1}

Rosen \cite{Ro54} established the $\lambda$-continued fractions
associated to the Hecke groups. Fixing $q\geq 3$ and 
$\lambda=2\cos\frac{\pi}{q}$ as always, a finite 
$\lambda$-continued fraction is a real number of the shape
\begin{eqnarray*}
[a_1;a_2,...,a_n]:=
a_1\lambda-\frac{1}{a_2\lambda-\frac{1}{a_3\lambda-...-
\frac{1}{a_n\lambda}}},
\end{eqnarray*}
where $a_1\in\Z$, $a_2,a_3,...,a_n\in\Z-\{0\}$. This extends
in the obvious way to infinite $\lambda$-continuous fractions.
The number above is
\begin{eqnarray*}
[a_1;a_2,...,an]=
\oline{A_1}^{a_1}\oline{V}\,\oline{A_1}^{a_2}\oline{V}...
\oline{V}\,\oline{A_1}^{a_n}\oline{V}(\infty).
\end{eqnarray*}
Obviously any cusp has a finite $\lambda$-continued fraction,
and any fixed point $r\in\R$ of an element 
$\oline{B}\in G_q$ has a periodic $\lambda$-continued fraction.

Rosen considered the {\it nearest integer algorithm}
for constructing a $\lambda$-continued fraction for any
real number $r$. We prefer to follow \cite{LL16}
and write it as {\it pseudo Euclidean algorithm}:
Write $r=\frac{r_0}{r_1}$ with $r_1\neq 0$ and 
construct a sequence $(r_0,r_1,r_2,...)$ of real numbers 
which either ends with some $r_{n+1}=0$ or never ends, 
in the following way. Given $r_i$ and $r_{i+1}\neq 0$
define $r_{i+2}\in\R$ and $a_{i+1}\in\Z$ uniquely by
\begin{eqnarray*}
r_i=a_{i+1}\lambda r_{i+1}-r_{i+2},\\
-|r_{i+1}|\frac{\lambda}{2}\leq r_{i+2}
<|r_{i+1}|\frac{\lambda}{2}.
\end{eqnarray*}
If $r_{i+2}\neq 0$ then
\begin{eqnarray*}
\frac{r_i}{r_{i+1}}=a_{i+1}\lambda - \frac{1}{r_{i+1}/r_{i+2}}.
\end{eqnarray*}
Observe $\frac{\lambda}{2}<1$, so 
$|r_{i+2}|\leq \frac{\lambda}{2}|r_{i+1}|<|r_{i+1}|$,
so the sequence $(|r_i|)_{i\geq 0}$ decreases exponentially
if it does not stop with a value 0. 

Rosen showed in \cite{Ro54}:

(1) The nearest integer algorithm gives for each cusp
a finite $\lambda$-continued fraction, and 

(2) The nearest integer algorithm gives for any real number 
which is not a cusp a convergent $\lambda$-continued fraction.

Part (1) follows now also easily with Theorem \ref{t1.1}:
Write a cusp $r$ as $r=\frac{r_0}{r_1}$ with 
$r_0+r_1\zeta^{q-1}\in\Delta^{(1)}_q$. Then also
$r_i+r_{i+1}\zeta^{q-1}\in\Delta^{(1)}_q$ as long as 
$r_{i+1}$ is defined.
Because $\Delta^{(1)}_q$ is discrete and especially does not accumulate
at 0, the sequence $(r_i)_{i\geq 0}$ stops with a value 0.

As said above, it is obvious that any $r\in\R$ which is a 
fixed point of a hyperbolic element of $G_q$ has a periodic
$\lambda$-continued fraction. But it is not clear whether
the nearest integer algorithm leads for such an $r$ to
a periodic $\lambda$-continued fraction.

A good point in the pseudo Euclidean algorithm above of
\cite{LL16} is the following: For each 
$(r_0,r_1)\in\Z[\lambda]$ it leads in a bounded number of steps
to a decision, whether $r_0+r_1\zeta^{q-1}\in\Delta^{(1)}_q$
or not. Though it does not lead in a bounded number of steps
to a decision whether $\frac{r_0}{r_1}\in G_q(\infty)$. 
Lang and Lang use this to characterize
the elements of $G_q^{mat}$: Consider $a,b,c,d\in\Z[\lambda]$.
Then \cite[Proposition 3.7]{LL16}
\begin{eqnarray*}\label{10.1} 
\begin{pmatrix}a&b\\c&d\end{pmatrix}\in G_q^{mat}&\iff&
a+c\zeta^{q-1}\in\Delta^{(1)}_q,b+d\zeta^{q-1}\in\Delta^{(1)}_q\\
&&\textup{and } ad-bc=1.
\end{eqnarray*}
This result seems to supersede \cite{Ro86} where Rosen 
considered the case $q=5$. 
The proof uses besides other lemmas Lemma 3.3 in \cite{LL16},
which says that $a+c\zeta^{q-1}\in\Delta^{(1)}_q$ satisfies
$|a|\geq 1$ and $|c|\geq 1$. This is part of our Theorem 
\ref{t1.1}. Also the construction of the pairs $(a,c)$ with 
$a+c\zeta^{q-1}\in\Delta^{(1)}_q$ in \cite[section 3.2]{LL16}
is related to our construction of $\Delta^{(1)}_q$ in 
section \ref{s5}. 

Rosen's paper \cite{Ro54} sets the stage for questions about
$G_q$-orbits in $\Q(\lambda)\cup\{\infty\}$.

\subsection{Hyperbolic fixed points}\label{s10.2}

Schmidt and Sheingorn give in the paper \cite{SS95}
in part 1.2 of the introduction a detailed survey on
results which connect the Hecke groups $G_q$ with geometry.

In the cases $q=4$ and $q=6$ Schmidt and Sheingorn 
\cite[Theorem 1]{SS95} show that all elements in
$\Q(\lambda)$ which are not cusps, so all elements of
$\Q(\lambda)-\lambda\Q$ are fixed points of hyperbolic elements. 
Theorem 1 in \cite{SS95} controls also all fixed points
in $\R-\Q(\lambda)$ of hyperbolic elements of $G_q$. 
It also characterizes all real numbers by $\lambda$-continued
fractions. But they do not use Rosen's nearest integer
algorithm, but a variant. With this variant precisely
the fixed points in $\R$ of parabolic or hyperbolic
elements of $G_q$ have finite or periodic 
$\lambda$-continued fractions. Though a slightly unpleasant
point is that with their variant a cusp can have a 
periodic $\lambda$-continued fraction.

\cite{RT01} and \cite{HMTY08} give interesting examples,
conjectures and results on orbits of $G_q$ in $\Q(\lambda)$.
For even $q$ \cite[Corollary 3]{HMTY08} says that $G_q$ 
has infinitely many orbits in $\Q(\lambda)\cup\{\infty\}$. 

For odd $q$ \cite[Theorem 9]{HMTY08} says that $G_q$ has
at least $(2^{\varphi(q)/2}+q+1)/2$ orbits in 
$\Q(\lambda)\cup\{\infty\}$. Table 1 in \cite{HMTY08} 
makes this lower bound for $q\in\{3,5,...,,37\}$ explicit.
Unfortunately, it differs from a similar table in \cite{BR73}.
Remark 2 in \cite{HMTY08} points at a possible mistake
in \cite{BR73}. 

On the other hand, \cite{BR73} as well as
\cite{HMTY08} are not careful with the assumption whether
$\Z[\lambda^2]$ is a principal ideal domain or not.
Lemma 6 in \cite{HMTY08} needs this assumption as it applies
the chinese remainder theorem, but does not formulate it.
Theorem 9 builds on Lemma 6.
Also the example at the bottom of page 84 in \cite{BR73}
needs this assumption, but does not make it.

How these orbits in \cite{HMTY08} look like, i.e. whether 
some of them consist of fixed points of hyperbolic elements, 
is not clear.
The following result is of interest in the context of
fixed points of hyperbolic elements.

\begin{theorem}\label{t10.1}\cite[Theorem 10.3.5]{Be83}
Each non elementary Fuchsian group contains infinitely many
conjugacy classes of maximal cyclic groups generated by
hyperbolic elements.
\end{theorem}

\cite[Theorem 3]{RT01} gives for $q=7$ many elements in 
$\Q[\lambda]$ where the $\lambda$-continued fraction from the 
nearest integer algorithm is periodic. Their periods
are all the same. Therefore they are all fixed points of 
hyperbolic elements, and these hyperbolic elements are conjugate.

\cite[4.1]{HMTY08} gives for $q=9$ four families of elements
where the $\lambda$-continued fractions from the nearest
integer algorithm are periodic, and the members of one family
have the same periods. Therefore they are all fixed points
of hyperbolic elements, and these hyperbolic elements
fall into four conjugacy classes.

These and other experiments lead them to the following questions:

(1) In the case $q=7$, does $\Q(\lambda)\cup\{\infty\}$ consist
only of two $G_q$-orbits, the orbit of cusps and one orbit of
hyperbolic fixed points?

(2) In the case $q=9$, does $\Q(\lambda)\cup\{\infty\}$ consist
only of five $G_q$-orbits, the orbit of cusps and four orbits of
hyperbolic fixed points?

(3) In the cases $q\in\{11,13,...,29\}$, does $\Q(\lambda)$
contain any orbits of fixed points of hyperbolic elements?
\cite{HMTY08} did not find any. (Of course, as $G_q$ contains
hyperbolic elements, $\R$ contains $G_q$-orbits of fixed points
of hyperbolic elements.)

\subsection{One more result for $q=5$}\label{s10.5}

McMullen \cite{Mc22}
studied the case $q=5$, coming from Veech surfaces
and Hilbert modular surfaces. Besides others, he has 
the following result: Let 
$a=\alpha_1+\alpha_2\lambda\in\Z[\lambda]$
with $\alpha_1,\alpha_2\in\Z$ be an entry of a matrix
in $G_5^{mat}$. Then 
$$\alpha_1\alpha_2\geq 0.$$

\subsection{Figures for $\Delta^{(1)}_q$ in the cases
$q\in\{8;9\}$}\label{s10.6}

\begin{figure}
\includegraphics[width=1.0\textwidth]{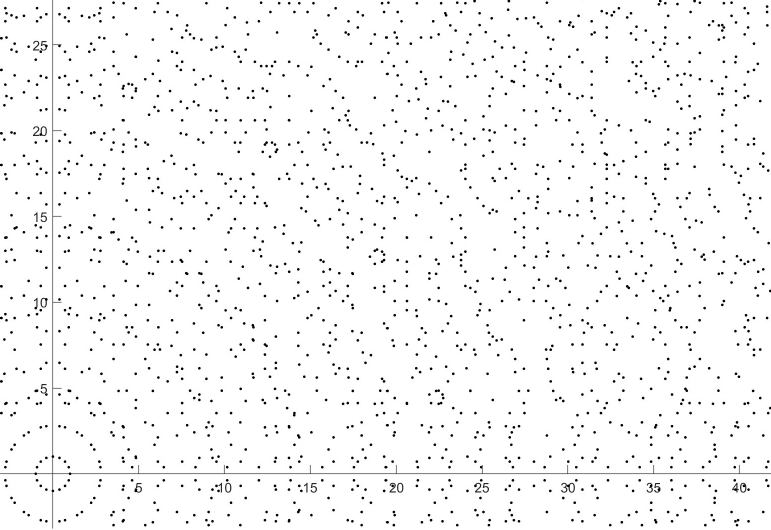}
\caption[Figure 10.1]{Part of the set $\Delta^{(1)}_8$}
\label{Fig:10.1}
\vspace*{0.5cm}
\end{figure}

\begin{figure}
\includegraphics[width=1.0\textwidth]{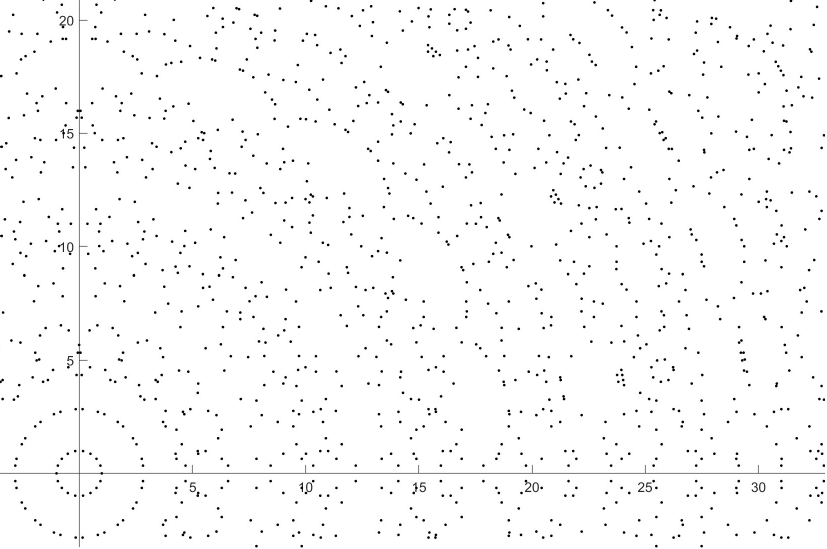}
\caption[Figure 10.2]{Part of the set $\Delta^{(1)}_9$}
\label{Fig:10.2}
\end{figure}

We finish the paper with the two Figures \ref{Fig:10.1} 
and \ref{Fig:10.2}. They show a part of the
set $\Delta^{(1)}_q$ in the cases $q=8$ and $q=9$. 

\begin{thebibliography}{AAA9}
\bibitem[AS09]{AS09} P. Arnoux, T.A. Schmidt: \quad
Veech surfaces with non-periodic directions in the trace field.
J. Mod. Dyn. {\bf 3.4} (2009), 611--629.
\bibitem[Be83]{Be83} A.F. Beardon: \quad 
The geometry of discrete groups.
Graduate Texts in Mathematics 91, Springer, 1983.
\bibitem[Bo73]{Bo73} W. Borho: \quad 
Kettenbrüche im Galoisfeld.
Abh. Math. Sem. Univ. Hamburg {\bf 39} (1973), 76--82.
\bibitem[BR73]{BR73} W. Borho, G. Rosenberger: \quad 
Eine Bemerkung zur Hecke-Gruppe $G(\lambda)$.
Abh. Math. Sem. Univ. Hamburg {\bf 39} (1973), 83--87. 
\bibitem[vdG88]{vdG88} G. van der Geer: \quad 
Hilbert modular surfaces. Springer, 1988.
\bibitem[HMTY08]{HMTY08} E. Hanson, A. Merberg, C. Towse, E. Yudovina: \quad
Generalized continued fractions and orbits under the action of Hecke triangle groups.
Acta Arithmetica {\bf 134.4} (2008), 337--348.
\bibitem[He36]{He36} E. Hecke: \quad
\"Uber die Bestimmung Dirichletscher Reihen durch ihre 
Funktionalgleichungen. Math. Ann. {\bf 112} (1936), 
664--699.
\bibitem[HL24]{HL24} C. Hertling, K. Larabi: \quad
Unimodular bilinear lattices, automorphism groups,
vanishing cycles, monodromy groups, distinguished bases,
braid group actions and moduli spaces from upper triangular
matrices.
Preprint arXiv:2412.16570, 396 pages. 
\bibitem[Hu90]{Hu90} J.E. Humphreys: \quad{}
Reflection groups and Coxeter groups. 
Cambridge Studies in Advanced Mathematics \textbf{29}.
Cambridge University Press, 1990. Duke Math. J. {\bf 123.1} (2004), 49--69.
\bibitem[Ko97]{Ko97} H. Koch: \quad 
Zahlentheorie. Algebraische Zahlen
und Funktionen. Vieweg, 1997.
\bibitem[LL16]{LL16} Ch.L. Lang, M.L. Lang: \quad
Arithmetic and geometry of the Hecke groups. 
J. Algebra {\bf 460} (2016), 392--417.
\bibitem[Le66]{Le66} J. Lehner: \quad 
A short course in automorphic functions.
Holt, Rinehart and Winston, New York, 1966.
\bibitem[Le67]{Le67} A. Leutbecher: \quad 
Über die Heckeschen Gruppen $G(\lambda)$.
Abh. Math. Sem. Univ. Hamburg {\bf 31} (1967), 119--205.
\bibitem[Le74]{Le74} A. Leutbecher: \quad 
Über die Heckeschen Gruppen $G(\lambda)$ II.
Math. Ann. {\bf 211} (1974), 63--84.
\bibitem[vdL82]{vdL82} F.J. van der Linden: \quad 
Class number computations of real abelian number fields.
Comp. Math {\bf 39.160} (1982), 693--707.
\bibitem[Mc22]{Mc22} C. McMullen: \quad
Billiards, heights, and the arithmetic of non-arithmetic groups.
Invent. Math. {\bf 228} (2022), 1309--1351.
\bibitem[Ro54]{Ro54} D. Rosen: \quad 
A class of continued fractions associated with certain properly
discontinuous groups.
Duke Math. J. {\bf 21} (1954), 549--564.
\bibitem[Ro63]{Ro63} D. Rosen: \quad 
An arithmetic characterization of the parabolic points of
$G(2\cos\frac{\pi}{5})$.
Glasgow Math.. Assoc {\bf 6} (1963), 88--96.
\bibitem[Ro86]{Ro86} D. Rosen: \quad 
The substitutions of the Hecke group $\Gamma(2\cos\frac{\pi}{5})$.
Arch. Math. {\bf 46} (1986), 533--538.
\bibitem[RT01]{RT01} D. Rosen, C. Towse: \quad 
Continued fraction representations of units associated with 
certain Hecke groups.
Arch. Math {\bf 77} (2001), 294--302.
\bibitem[SS95]{SS95} T.A. Schmidt, M. Sheingorn: \quad 
Length spectra of the Hecke groups.
Math. Z. {\bf 220} (1995), 369--397.
\bibitem[Sch02]{Sch02} R. Schoof: \quad
Class numbers of real cyclotomic fields of prime conductor.
Math. Comp {\bf 72.242} (2002), 913--937. 
\bibitem[Se85]{Se85} F. Seibold: \quad
Zahlentheoretische Eigenschaften der Heckeschen Gruppen $G(\lambda)$
und verwandter Transformationsgruppen.
Inauguraldissertation Techn. Univ. M\"unchen 1985 
(Betreuer A. Leutbecher).
\bibitem[Vi71]{Vi71} \`E.B. Vinberg: \quad
Discrete linear groups generated by reflections.
Math. USSR Izvestija {\bf 5.5} (1971), 1083--1119.
\bibitem[Wa82]{Wa82} L.C. Washington: \quad 
Introduction to cyclotomic fields.
Graduate texts in mathematics {\bf 83}, Springer 1982.
\bibitem[Wo77]{Wo77} J. Wolfart: \quad
Eine Bemerkung über Heckes Modulgruppen.
Arch. Math. {\bf 29} (1977), 72--77. 
\end{thebibliography}
\end{document}